\documentclass{article}  
\usepackage{amssymb}              
\usepackage{amsthm}
\usepackage{amsmath}
\usepackage{eucal}
\usepackage{verbatim}
\usepackage{graphicx}
\usepackage{multirow}
\usepackage{fancyhdr}
\usepackage{color}
\usepackage{enumerate}
\usepackage{calrsfs}
\usepackage{cancel}
\usepackage{fullpage}
\usepackage{amsmath}

\newtheorem{theorem}{Theorem}
\newtheorem{lemma}{Lemma}
\newtheorem{definition}{Definition}
\newtheorem{remark}{Remark}
\newtheorem{proposition}{Proposition}
\newtheorem{corollary}{Corollary}
\newtheorem{example}{Example}

\makeatletter
\newcommand{\leqnomode}{\tagsleft@true}
\newcommand{\reqnomode}{\tagsleft@false}
\makeatother

\def\({\begin{eqnarray}}
\def\){\end{eqnarray}}
\def\[{\begin{eqnarray*}}
\def\]{\end{eqnarray*}}
\def\part#1#2{\frac{\partial #1}{\partial #2}}
\def\parts#1#2#3{\frac{\partial^2 #1}{\partial #2 \partial #3}}

\def\R{\mathbb{R}}
\def\N{\mathbb{N}}

\def\eps{\varepsilon}

\def\Norm#1{\left\| #1 \right\|}

\def\sign{\mbox{sign}}

\def\Lap{\mathcal{L}}

\newcommand\Vset{V}
\newcommand\Eset{E}
\newcommand\Cset{\mathcal{C}}

\def\Ekin{\E_{\mathrm{kin}}}
\def\Emet{\E_{\mathrm{met}}}

\def\E{\mathcal{E}}
\def\F{\mathcal{F}}

\def\f{\mathfrak{f}}
\def\chee{\mathfrak{h}}
\def\minL{\ell}

\title{Robust network formation with biological applications}

\date{}

\begin{document}
\pagenumbering{arabic}

\author{Jan Haskovec\footnote{Mathematical and Computer Sciences
            and Engineering Division,
         King Abdullah University of Science and Technology,
         Thuwal 23955-6900, Kingdom of Saudi Arabia;
         {\it jan.haskovec@kaust.edu.sa}}\ \ and Jan Vyb\'\i ral\footnote{Department of Mathematics, Faculty of Nuclear Sciences and Physical Engineering,
        Czech Technical University, Trojanova 12, 12000 Praha, Czech Republic;
         {\it jan.vybiral@fjfi.cvut.cz}.
         The work of this author has been supported by the grant P202/23/04720S of the Grant Agency of the Czech Republic.}}

\maketitle
\vspace{2mm}

\textbf{Abstract.} We provide new results on the structure of optimal transportation networks obtained as minimizers of an energy cost functional consisting of a kinetic (pumping) and material (metabolic) cost terms, constrained by a local mass conservation law.
In particular, we prove that every tree (i.e., graph without loops) represents a local minimizer of the energy with concave metabolic cost. For the linear metabolic cost, we prove that the set of minimizers contains a loop-free structure.
Moreover, we enrich the energy functional such that it accounts also for robustness of the network, measured in terms of the Fiedler number of the graph with edge weights given by their conductivities.
We examine fundamental properties of the modified functional, in particular, its convexity and differentiability. We provide analytical insights into the new model by considering two simple examples. Subsequently, we employ the projected subgradient method to find global minimizers of the modified functional numerically.
We then present two numerical examples, illustrating how the optimal graph's structure and energy expenditure depend on the required robustness of the network.
\vspace{2mm}

\vspace{2mm}

\vspace{3mm}

\section{Introduction}
In this paper we focus on the discrete graph model introduced by Hu and Cai \cite{Hu-Cai} and further studied in \cite{bookchapter, Hu-Cai-19, BHMR, HKM}
which describes optimal transportation networks in primarily biological context.
Typical examples are leaf venation in plants, mammalian circulatory systems that convey nutrients to the
body through blood circulation, or neural networks that transport electric charge.
Understanding the properties of optimal transportation networks and mechanisms of their development, function and adaptation
has been subject of an active field of research, see, e.g., \cite{Bebber, Cantarella, Runions, Yancopoulos}.

Mathematical modeling of transportation networks is traditionally based on the frameworks of mathematical graph theory and discrete energy optimization.
The model \cite{Hu-Cai} falls into this category, being formulated as an energy functional depending on edge conductivities of a given undirected discrete graph,
constrained by a local mass conservation law. The mass conservation law is imposed in terms of a linear system of equations for the material pressures known as Kirchhoff law.
The energy consists of a pumping term, describing the kinetic energy of the material flow through the network,
and a metabolic term, which reflects the biological motivation of the model.
The metabolic term is a function of the edge conductivities and is
assumed to be of power-law form, with exponent $\gamma>0$. For biological applications,
one usually assumes that $\gamma\in [1/2,1]$. We note that for $\gamma<1$ the energy functional is nonconvex,
while for $\gamma=1$ it is convex and for $\gamma>1$ it is strictly convex.

The first goal of this paper is to provide two new results regarding the (local) minimizers of the energy.
In particular, it has been shown \cite{BHMR} that for $\gamma <1$, each minimizer of the energy is a loop-free graph, i.e., an undirected graph where each pair of nodes is connected by at most one path.
Here we show the complementary claim, namely, that if $\gamma<1$, every admissible loop-free graph represents a local minimizer
of the energy. The term ``admissible" refers to a graph for which the Kirchhoff law is solvable
when all edges with zero conductivity
are treated as void (nonexistent).
Moreover, for $\gamma=1$ we show that every (local) minimizer is either loop-free,
or there exist a matrix of conductivities with the same value of the energy, which represents a loop-free graph.

The second goal of this paper is to incorporate the concept of robustness into the model \cite{Hu-Cai}.
In particular, the fact that, for biological applications where $\gamma\leq 1$, the energy minimizers are loop-free,
does not correspond to the transportation networks observed in most real organisms,
where typically (many) loops are present. For a striking example of a highly redundant pattern in leaf venation
we refer to \cite[Figure 1]{Laguna-Bohn-Jagla}.
There the authors make the point that the high redundancy of paths from the leaf base to any point on the leaf
surface might be very advantageous with regard to local damages,
i.e., they promote robustness of the transportation network
against removal of edges.
We therefore propose to incorporate the robustness aspect into the discrete model \cite{Hu-Cai}
by extending the energy functional by a term that measures the connectivity of the network.

A generic measure of connectivity of a discrete graph is the Cheeger constant \cite{Mohar}, also called isoperimetric number.
Loosely speaking, the Cheeger constant is a numerical measure of whether or not a graph has a ``bottleneck",
in the sense of its resilience against disturbance of connectivity through removal of edges.
The Cheeger constant is strictly positive if and only if the underlying (undirected) graph is connected.
Graphs with small (but positive) Cheeger constant have a "bottleneck", in the sense that there are two "large" sets of vertices with "few" links (edges) between them.
On the other hand, the Cheeger constant is "large" if any possible division of the vertex set into two subsets has "many" links between those two subsets.
The fundamental problem is that the calculation of the isoperimetric number of graphs with multiple edges is NP hard \cite{Mohar}.
To circumvent this issue and to obtain a practically tractable optimization problem,
we propose to replace it by a related quantity, the algebraic connectivity, also called Fiedler number.
This is defined as the second smallest (i.e., for connected graph, the smallest nonzero) eigenvalue of the matrix Laplacian
of the connectivity matrix of the graph  \cite{Fiedler}. The Fiedler number is related to the Cheeger constant by the well known Cheeger inequalities \cite{Mohar},
which state that the isoperimetric number is bounded from below by one half of the Fiedler number.
Consequently, knowing that a given graph is ``well connected" due to sufficiently large value of its Fiedler number
guarantees robustness also in terms of its Cheeger constant. 

In this paper we take into account the connectivity of the graph
in terms of the Fiedler number calculated from the weighted adjacency matrix,
with weights equal to the edge conductivities.
This is motivated by the fact that, if a certain transport path
is severed due to removal of some edges, then the ``alternative" (backup) path should have
enough transport capacity, i.e., high enough conductivity.
The Fiedler number
is then a Lipschitz continuous function of the edge conductivities.
Moreover, as long as it represents a simple eigenvalue (i.e., of multiplicity one) of the matrix Laplacian,
it is differentiable with respect to the conductivities.

We thus formulate a new energy functional, where a multiple of the Fiedler number is subtracted
from the pumping (kinetic) and metabolic energy terms.
We then present two toy examples demonstrating how the modified energy functional
propagates connectivity of the optimal transportation network.
Moreover, since the Fiedler number is a concave function of the edge conductivities, the problem remains convex as long as $\gamma\geq 1$. The convexity and Lipschitz continuity of the functional facilitates the application of the projected subgradient method for its numerical minimization.
We then present results of two numerical experiments, further documenting how the connectivity of the graph and its energy expenditure
grows with increasing relative weight of the Fiedler number in the modified energy functional.

This paper is organized as follows. In Section \ref{sec:model} we describe the discrete network formation model of Hu and Cai \cite{Hu-Cai} and establish some of its fundamental properties (convexity, boundedness of the kinetic energy). We also introduce the matrix Laplacian and Fiedler number of the graph. In Section \ref{sec:structure} we provide results about the structure of optimal transportation networks for values of $\gamma<1$ and $\gamma=1$.
In Section \ref{sec:robustness} we extend the energy functional such that it also accounts for robustness of the transportation network in terms of the Fiedler number of the weighted graph. We then establish basic properties of the modified functional with $\gamma=1$, namely, boundedness from below, coercivity and convexity. We also provide insights into the new model by studying two toy examples. Finally, in Section \ref{sec:num} we implement the projected subgradient method to find global minimizers of the modified functional numerically. We present two examples where we demonstrate how the structure of the optimal graphs depends on the required robustness of the network.

\section{The model}\label{sec:model}
The discrete network formation model introduced by Hu and Cai \cite{Hu-Cai} is posed on a prescribed undirected connected graph
$G=(\Vset,\Eset)$, consisting of a finite set of vertices $\Vset$, also called nodes, and a finite set of edges $\Eset$.
The number of vertices shall be denoted by  $|\Vset|$ in the sequel.
Any pair of vertices is connected by at most one edge and no vertex is not connected to itself by an edge (i.e., the graph does not contain ``self-loops").
The edge  between vertices $i\in\Vset$ and $j\in\Vset$ is denoted by $(i,j)\in \Eset$. Since the graph is undirected we refer by $(i,j)$ and $(j,i)$ to the same edge.
For each edge $(i,j)\in\Eset$ of the graph $G$ we prescribe its length $L_{ij}=L_{ji}>0$.

Throughout the paper we assume that the graph $G$ is connected, i.e., that every pair of vertices in $V$ is connected by an undirected path (sequence of edges) in $E$.

We treat the graph as a transportation structure and
in each vertex $i\in\Vset$ we have the pressure $P_i\in\R$ of the medium flowing through the network.
The oriented flux (flow rate) from vertex $i\in\Vset$ to $j\in\Vset$ is denoted by $Q_{ij}$; obviously, we have the antisymmetry $Q_{ij}=-Q_{ji}$.
For biological networks, the Reynolds number of the flow is typically small and the flow is predominantly in the laminar (Poiseuille) regime.
Then the flow rate between vertices $i\in\Vset$ and $j\in\Vset$ along the edge $(i,j)\in\Eset$ is proportional to the conductance $C_{ij}\geq 0$ of the edge and the pressure drop $P_i-P_j$,
\begin{align}\label{eq:flowrate}
   Q_{ij} := C_{ij}\frac{P_i-P_j}{L_{ij}}\qquad\text{for all~}(i,j)\in \Eset.
\end{align}
The local mass conservation in each vertex is expressed in terms of the Kirchhoff law
\begin{align}\label{eq:kirchhoff}
   \sum_{j\in N(i)} C_{ij}\frac{ P_i-P_j}{L_{ij}}=S_i\qquad \text{for all~}i\in \Vset.
\end{align}
Here $N(i)$ denotes the set of vertices connected to $i\in\Vset$ through an edge, i.e.,
$N(i) := \{ j\in\Vset;\, (i,j)\in\Eset \}.$
Moreover, $S=(S_i)_{i\in\Vset}$ is the prescribed strength of the flow source ($S_i>0$) or sink ($S_i<0$) at vertex $i$.
Given the matrix of nonnegative conductivities $C=(C_{ij})_{(i,j)\in\Eset}$ and the prescribed positive lengths $L=(L_{ij})_{(i,j)\in\Eset}$, the Kirchhoff law \eqref{eq:kirchhoff} is a linear system of equations for the vector of pressures $P=(P_i)_{i\in\Vset}$.
Clearly, a necessary condition for the solvability of \eqref{eq:kirchhoff} is the global mass conservation
$S\in\R^{|V|}_0$, where here and in the sequel we denote
\begin{equation}\label{eq:RV0}
\R^{|V|}_0 := \left\{ v\in\R^{|V|}; \; \sum_{i=1}^{|V|} v_i = 0 \right\}.
\end{equation}

\begin{remark}\label{rem:Kirchhoff}
If $C$ is the adjacency matrix of a weighted connected graph, then
the solution $P$ of \eqref{eq:kirchhoff} always exists and is unique up to an additive constant.
Although $E$ itself is connected, $C$ might be supported on a proper subset of $E$ and may therefore correspond to a 
disconnected graph;
we then shortly say that $C$ is disconnected.
In this case, the solution of \eqref{eq:kirchhoff} may not exist or, if it exists,
it is not unique (up to an additive constant),
as demonstrated in Example \ref{ex:P} in the Appendix (Section \ref{sec:app}).
In case when $C$ is disconnected and
\eqref{eq:kirchhoff} is solvable, we denote the connected components of
$C$ by $\Gamma_1,\dots,\Gamma_n$ and observe that $\sum_{j\in \Gamma_k}S_j=0$ for every $1\le k\le n$.
Then $P$ is uniquely defined up to $n$ additive constants, one for every component $\Gamma_k$.
Therefore, the matrix of fluxes $Q$ is still uniquely defined by \eqref{eq:flowrate}.
On the other hand, if \eqref{eq:kirchhoff} is not solvable, then $C$ is disconnected.
\end{remark}


The conductivities $C_{ij} = C_{ji} \geq 0$ of the edges are subject to an energy optimization process.
Hu and Cai \cite{Hu-Cai} propose an energy cost functional consisting of a pumping power term and a metabolic cost term.
According to the Joule's law, the power (kinetic energy) needed to pump material through an edge $(i,j)\in\Eset$ is proportional to the pressure drop $P_i-P_j$
and the flow rate $Q_{ij}$ along the edge, i.e.,
\begin{align*}
  (P_i-P_j) Q_{ij}=\frac{Q_{ij}^2}{C_{ij}}L_{ij}. 
\end{align*}
The metabolic cost of maintaining the edge is assumed to be proportional to its length $L_{ij}$ and a power of its conductivity $C_{ij}^{\gamma}$, with an exponent $\gamma>0$.
For instance, in blood vessels the metabolic cost is proportional to the cross-section area of the vessel \cite{Murray}. Modeling the blood flow by Hagen-Poiseuille's law,
the conductivity is proportional to the square of the cross-section area, implying $\gamma=1/2$ for blood vessel systems.
For models of leaf venation the material cost is proportional to the number of small tubes, which is
proportional to $C_{ij}$, and the metabolic cost is due to the effective loss of the photosynthetic power
at the area of the venation cells, which is proportional to $C_{ij}^{1/2}$. Consequently, the effective value
of $\gamma$ typically used in models of leaf venation lies between $1/2$ and $1$, see, e.g., \cite{Hu-Cai-19}.
The energy cost functional is thus given by
\begin{align}\label{eq:energy}
   {\E}[C] := \sum_{(i,j)\in\Eset}\left( \frac{Q_{ij}[C]^2}{C_{ij}}+\frac{\nu}{\gamma} C_{ij}^{\gamma}\right) L_{ij},
\end{align}
where $Q_{ij}=Q_{ij}[C]$ is given by \eqref{eq:flowrate} with pressures calculated from the Kirchhoff law \eqref{eq:kirchhoff},
and $\nu>0$ is the so-called metabolic coefficient.
Note that  every edge of the graph $G$ is counted exactly once in the above sum,
i.e., we identify each edge $(i,j)$ with $(j,i)$ and the energy can also be written as
\[
   {\E}[C] := \frac12 \sum_{i\in\Vset} \sum_{j\in\Vset} \left( \frac{Q_{ij}[C]^2}{C_{ij}}+\frac{\nu}{\gamma} C_{ij}^{\gamma}\right) L_{ij},
\]
where we set $L_{ij} := 0$ whenever $(i,j)\notin\Eset$.
If $C$ is such that the linear system \eqref{eq:kirchhoff} is not solvable, we formally set $\E[C]:=+\infty$.
In the sequel, we shall sometimes address the kinetic (pumping) and metabolic parts of the energy separately. For this purpose we denote
\(
   \label{eq:Ekinmet}
   \Ekin[C] := \sum_{(i,j)\in \Eset} \frac{Q_{ij}[C]^2}{C_{ij}} L_{ij}, 
   \qquad
   \Emet[C] := \frac{\nu}{\gamma} \sum_{(i,j)\in \Eset} C_{ij}^\gamma L_{ij}.
\)

In order to find the optimal transportation structure for a given vector of sources and sinks $S\in\R^{|V|}_0$, one needs to find the (global) minimum of the energy functional \eqref{eq:energy}, coupled to the Kirchhoff law \eqref{eq:flowrate}--\eqref{eq:kirchhoff}.
The energy functional is to be minimized over the convex set $\Cset$ of symmetric matrices with nonnegative elements,
\( \label{eq:setC}
    \Cset := \left\{ C\in\R^{|V|\times |V|}; \, C_{ij}=C_{ji} \geq 0 \mbox{ for all } i, j\in V, \mbox{ with } C_{ij}=C_{ij}=0 \mbox{ if } (i,j)\notin\Eset \right\}.
\)
If the optimal solution has $C_{ij}=C_{ij}=0$ for any $(i,j)\in\Eset$, then the corresponding edge $(i,j)$
is considered nonexistent. We collected some basic mathematical properties of the functional $\E=\E[C]$
in the Appendix (Section \ref{sec:app}).
An important information in the context of the optimization task \eqref{eq:flowrate}--\eqref{eq:energy} regards the convexity properties
of the problem.

\subsection{Convexity of the energy functional}

\begin{proposition}\label{lem:convex}
The energy functional \eqref{eq:energy}, constrained by the Kirchhoff law \eqref{eq:flowrate}--\eqref{eq:kirchhoff},
is strictly convex for $\gamma>1$ and convex for $\gamma=1$.
\end{proposition}

This result is not new, see for instance \cite{HMP22, HMP23} for an analogous result
for the continuum version of the problem, but for the discrete case it was never explicitly formulated or proven in any of the papers \cite{bookchapter, BHMR, Hu-Cai-19, HKM}.
It follows directly from the convexity of the pumping (kinetic) part of the energy
$\Ekin[C]$,
coupled to the Kirchhoff law \eqref{eq:flowrate}--\eqref{eq:kirchhoff}. If \eqref{eq:kirchhoff} is not solvable, we again set $\Ekin[C]:=+\infty.$
We provide the proof here for the sake of the reader.

\begin{lemma}\label{lem:kin_convex}
The pumping energy $\Ekin[C]$ defined in \eqref{eq:Ekinmet}, constrained by the Kirchhoff law \eqref{eq:flowrate}--\eqref{eq:kirchhoff},
is a convex functional on the set $\Cset$.
\end{lemma}

\begin{proof}
By Lemma \ref{lem:C1} of the Appendix, the set $\{C\in \Cset:\E(C)<+\infty\}$ is an open, convex subset of $\Cset$
and the same obviously holds true with $\E$ replaced by $\Ekin$. To show the convexity of $\Ekin$, it is therefore enough
to prove that the Hessian matrix of $\Ekin[C]$ is positive semidefinite for every $C\in \Cset$ with $\Ekin[C]<+\infty.$

We use \eqref{eq:flowrate} to express the kinetic energy as
\(   \label{eq:Ekin}
   \Ekin[C] = \sum_{(i,j)\in \Eset} C_{ij} \frac{(P_j-P_i)^2}{L_{ij}}.
\)
We note that, by Lemma \ref{lemma:pert} of the Appendix, $P$ restricted by the condition $\sum_{j\in V} P_j=0$ is a differentiable function of $C$.
The first-order derivative of $\Ekin[C]$ with respect to the element $C_{km}$ reads, see \cite[Lemma 2.1]{HKM},
\( \label{eq:derivative}
   \part{}{C_{km}} \Ekin[C] = 
      - \frac{(P_k - P_m)^2}{L_{km}}.
\)
Consequently, the second-order derivative with respect to the elements $C_{km}$ and $C_{\alpha\beta}$ reads
\[
   \parts{}{C_{km}}{C_{\alpha\beta}}  \Ekin[C] = - \frac{1}{L_{km}}  \part{}{C_{\alpha\beta}} (P_k - P_m)^2.
\]
We fix a vector $\varphi\in\R^{|V|}$, multiply the Kirchhoff law \eqref{eq:kirchhoff} by $\varphi_i$ and sum over $i\in \Vset$, using the standard symmetrization trick on the left-hand side,
\[
   \frac12 \sum_{i\in \Vset} \sum_{j\in \Vset} \frac{ C_{ij}}{L_{ij}} (P_i-P_j) (\varphi_i-\varphi_j) = \sum_{i\in \Vset} S_i \varphi_i.
\]
We take a derivative of the above identity with respect to $C_{km}$,
\[
   \frac{1}{L_{km}} (P_k-P_m)(\varphi_k-\varphi_m) + 
      \frac12 \sum_{i\in \Vset} \sum_{j\in \Vset} \frac{ C_{ij}}{L_{ij}}  (\varphi_i-\varphi_j) \part{}{C_{km}} (P_i-P_j) = 0,
\]
where we took into account the symmetry $C_{km} = C_{m k}$.
We now choose $\varphi := \part{}{C_{\alpha\beta}} P$, which gives
\[
   \frac{1}{L_{km}} (P_k-P_m) \part{}{C_{\alpha\beta}} (P_k-P_m) + 
      \frac12 \sum_{i\in \Vset} \sum_{j\in \Vset} \frac{ C_{ij}}{L_{ij}}  \part{}{C_{\alpha\beta}} (P_i-P_j) \part{}{C_{km}} (P_i-P_j) = 0.
\]
Consequently,
\[
    \frac{1}{L_{km}} \part{}{C_{\alpha\beta}} (P_k-P_m)^2 = - \sum_{i\in \Vset} \sum_{j\in \Vset} \frac{ C_{ij}}{L_{ij}}  \part{}{C_{\alpha\beta}} (P_i-P_j) \part{}{C_{km}} (P_i-P_j)
\]
and
\[
   \parts{}{C_{km}}{C_{\alpha\beta}} \Ekin[C] = \sum_{i\in \Vset} \sum_{j\in \Vset} \frac{ C_{ij}}{L_{ij}}  \part{}{C_{\alpha\beta}} (P_i-P_j) \part{}{C_{km}} (P_i-P_j).
\]
Now, fix any $\xi\in \R^{|V|\times|V|}$ and denote
\[
   \Xi_{ij} := \sum_{k\in\Vset} \sum_{m\in\Vset}  \xi_{km} \part{}{C_{km}} (P_i-P_j).
\]
Then we have
\[
   \sum_{k, m} \sum_{\alpha,\beta} \left( \parts{\Ekin[C] }{C_{km}}{C_{\alpha\beta}} \right) \xi_{km} \xi_{\alpha\beta}
      = \sum_{i\in \Vset} \sum_{j\in \Vset} \frac{ C_{ij}}{L_{ij}} \Xi_{ij}^2 \geq 0,
\]
where the nonnegativity follows from $C_{ij}\geq 0$.
We conclude that the Hessian matrix of $\Ekin[C]$ is positive semidefinite for every $C\in \Cset$ with $\Ekin[C]<+\infty$, and, therefore, $\Ekin$ is convex in $\Cset$.
\end{proof}

\subsection{Matrix Laplacian and Fiedler number}

\begin{definition} \label{def:Lapl}
For any $C\in\Cset$ we define the matrix Laplacian $\Lap[C]$ as
\[
   \Lap[C] := D - C,
\]
where $D$ is the diagonal matrix of row/column sums of $C$ (recall that $C\in\Cset$ is symmetric).
If no risk of confusion, we shall often write just $\Lap$ instead of $\Lap[C]$.
\end{definition}

With this definition, we can express the Kirchhoff law \eqref{eq:kirchhoff} in the form
\(  \label{eq:kLap}
   \Lap[\widetilde C] P = S,
\)
where we introduced the matrix $\widetilde C_{ij} := C_{ij}/L_{ij}$, with $\widetilde C_{ij} :=0$ if $L_{ij}=0$
(recall that $L_{ij}=0$ means that $(i,j)\notin E$ and, consequently, $C_{ij}=0$).
Moreover, we denoted $P=(P_i)_{i\in\Vset}$ the vector of pressures
and $S=(S_i)_{i\in\Vset}$ the vector of sources and sinks.

We also note that if $P$ is a solution of \eqref{eq:kLap}, then an easy calculation gives the following expression for the kinetic energy \eqref{eq:Ekinmet},
\(  \label{EkinLap}
   \Ekin[C] = P^T \Lap[\widetilde C] P.
\)

\begin{definition} \label{def:Fied}
The Fiedler number $\f[C]$ of a matrix $C\in\Cset$ is the second smallest eigenvalue of the matrix Laplacian $\Lap[C]$,
where multiple eigenvalues are counted separately \cite{Fiedler}.
\end{definition}

\subsection{An upper bound on the kinetic energy}

\begin{lemma}\label{lem:ekinbound}
Let $0 \neq S \in R_0^{|V|}$.
For any $C\in\Cset$ we have
\(   \label{eq:ekinbound}
   \Ekin[C] \leq \frac{\Norm{S}^2}{\f[\widetilde C]},
\)
where $\Ekin[C]$ is defined in \eqref{eq:Ekinmet} and
$\widetilde C_{ij} := C_{ij}/L_{ij}$, with $\widetilde C_{ij} :=0$ if $L_{ij}=0$.
Here and in the sequel, $\Norm{S}$ denotes the $\ell^2$-norm of the vector $S$.
\end{lemma}

\begin{proof}
For the case when $C\in\Cset$ is such that the Kirchhoff law \eqref{eq:kLap} is not solvable, we have set $\Ekin[C] = +\infty$.
But then the graph composed of edges $(i,j)$ with $C_{ij}>0$ is not connected, see Remark \ref{rem:Kirchhoff}, thus $\f[C]=0$ and \eqref{eq:ekinbound} holds.

Let us now fix $C\in\Cset$ such that the Kirchhoff law \eqref{eq:kLap} admits a solution $P\in\R^{|V|}$. Since $P$ can be shifted by an additive constant,
we choose $P\in\R^{|V|}_0$. Therefore, denoting $\mathbf{1} := (1, \ldots, 1)\in\R^{|V|}$, we have
$P^T \mathbf{1} = 0$.
The Courant-Fisher (minimax) theorem \cite{Parlett} gives
\(  \label{Courant}
   \f[\widetilde C] = \min \left\{ \frac{v^T 
   \Lap[\widetilde C]v}{\Norm{v}^2}; \, v\in\R^{|V|}, \, v^T \mathbf{1} = 0 \right\}.
\)
Consequently,
\(  \label{eq:PP}
   \f[\widetilde C] \Norm{P}^2 \leq
   P^T 
   \Lap[\widetilde C] P.
\)
Taking a scalar product of \eqref{eq:kLap} with $P$ and applying the Cauchy-Schwartz and Young inequalities yields
\[
   P^T 
   \Lap[\widetilde C] P = P^T S \leq \Norm{S} \Norm{P} \leq \frac{\Norm{S}^2}{2 \f[\widetilde C]} + \frac{\Norm{P}^2}{2} \f[\widetilde C].
\]
Combining this with \eqref{eq:PP}, we obtain
\[
   P^T \Lap[\widetilde C] P \leq \frac{\Norm{S}^2}{\f[\widetilde C]}
\]
and we conclude by using \eqref{EkinLap}.
\end{proof}

\section{Structure of optimal transportation networks for $\gamma\leq 1$}\label{sec:structure}

In this section we shall provide results about the structure of optimal transportation networks for $\gamma<1$ and $\gamma=1$. Let us recall that we assume the graph $G=(\Vset, \Eset)$ to be connected.

\subsection{Trees are local minimizers for $\gamma<1$}\label{subsec:trees}

In \cite[Theorem 2.1]{BHMR} the authors proved that every global minimizer of the energy functional \eqref{eq:energy} with $\gamma<1$, constrained by the Kirchhoff law \eqref{eq:flowrate}--\eqref{eq:kirchhoff}, is loop-free (again, we remind that edges with zero conductivities are treated as nonexistent).
By a very slight modification of the proof, one can extend the claim to local minimizers.
Here we prove the opposite claim, namely, that for $\gamma<1$ and every spanning tree of the connected graph $(\Vset, \Eset)$ we can find a set of conductivities $C\in\Cset$ supported on the spanning tree which is a local minimizer of the constrained energy.

\begin{theorem}\label{thm:new:1}
Let $\gamma<1$ and let $(V,\tilde E)$, $\tilde E \subseteq E$, be a spanning tree of the connected graph $(V,E)$. Then there exists a set of conductivities $\tilde C\in \Cset$, supported on $\tilde E$, which is a local minimizer of the energy \eqref{eq:energy} constrained by the Kirchhoff law \eqref{eq:flowrate}--\eqref{eq:kirchhoff}, on the set $\Cset$.
\end{theorem}

\begin{proof}
We observe that since $(V,\tilde E)$ is a tree, the fluxes $\tilde Q_{ij}$ for $(i,j)\in\tilde E$ are uniquely determined by $S$.
Indeed, for any edge $(i,j)\in\tilde E$, there is a unique split of the vertices into two disjoint sets $V^{(i)}$, $V^{(j)}$ such that nodes from $V^{(i)}$ are connected to $i$ by paths in $\tilde E \setminus \{(i,j)\}$, and analogously for nodes in $V^{(j)}$.
By the local mass conservation, the flux through $(i,j)$ is then given by
\( \label{eq:tildeQij}
   \tilde Q_{ij} = \sum_{k\in V^{(i)}} S_k - \sum_{k\in V^{(j)}} S_k.
\)
Once we identified the fluxes by the above prescription, we construct the conductivities $\tilde C\in\Cset$ as follows: We set $\tilde C_{ij}:=(\tilde Q_{ij}^2/\nu)^{1/(\gamma +1)}$ for all $(i,j)\in\tilde\Eset$, and $\tilde C_{ij}:=0$ for $(i,j)\in E\setminus\tilde\Eset$. Obviously, $\tilde C$ is supported on $\tilde\Eset$.
We now claim that $\tilde C$ represents a local minimizer of the energy \eqref{eq:energy} constrained by the Kirchhoff law \eqref{eq:flowrate}--\eqref{eq:kirchhoff} on the set $\Cset$.

Let us first observe that since $(\Vset,\tilde\Eset)$ is connected, the Kirchhoff law \eqref{eq:kirchhoff} admits a solution $\tilde P\in\R^{|V|}$. Clearly, the fluxes $\tilde Q_{ij}$ constructed in \eqref{eq:tildeQij} verify the relation \eqref{eq:flowrate}. Moreover, we recall the formula \cite[Lemma 2.1]{HKM} for the derivative
\begin{equation}\label{eq:diffE}
   \frac{\partial {\mathcal E}[\tilde C]}{\partial \tilde C_{ij}}=-\left(\frac{ \tilde Q_{ij}^2}{\tilde C_{ij}^2}-\nu \tilde C_{ij}^{\gamma-1}\right)L_{ij}.
\end{equation}
Then, by construction, we have
\[
   \frac{\partial {\mathcal E}[\tilde C]}{\partial \tilde C_{ij}}=0 \qquad
   \mbox{for all } (i,j)\in\tilde\Eset \mbox{ with } \tilde C_{ij}>0.
\]
Moreover, using \eqref{eq:flowrate} we recast \eqref{eq:diffE} as
\[
   \frac{\partial \E[\tilde C]}{\tilde C_{ij}}=-\frac{(\tilde P_j-\tilde P_i)^2}{L_{ij}}+\nu \tilde C_{ij}^{\gamma-1}L_{ij},
\]
which gives
\[
   \frac{\partial {\mathcal E}[\tilde C]}{\partial \tilde C_{ij}}=+\infty \qquad
   \mbox{for all } (i,j)\in\Eset \mbox{ with } \tilde C_{ij}=0,
\]
where the partial derivative is taken from the right, i.e., as $\tilde C_{ij}\to 0+$.
Since $\tilde C\in\Cset$ is unique element of $\Cset$ with these properties, we conclude that it is
a local minimizer of $\E=\E[C]$ on $\Cset$.
\end{proof}

\begin{remark}
Combining the claims of Theorem \ref{thm:new:1} and \cite[Theorem 2.1]{BHMR} implies that the set of local minimizers of 
the energy functional \eqref{eq:energy} with $\gamma<1$ is the set of spanning trees of the graph $(\Vset,\Eset)$. Consequently, a possible way to find the global minimizer is to search the set of spanning trees and identify the one with the smallest value of the energy $\E=\E[C]$. The value of the energy for any given tree is calculated using the procedure described in the proof of Theorem \ref{thm:new:1}. Clearly, this approach quickly turns computationally infeasible with growing size of the graph. For instance, the number of spanning trees for a complete graph with $|\Vset|$ nodes is $|\Vset|^{|\Vset|-2}$ by the Cayley's formula. For planar graphs, the number of possible spanning trees grows exponentially, see, e.g., \cite{Buchin-Schulz}.
\end{remark}

\subsection{Structure of the set of minimizers for $\gamma=1$}

In this section we consider the energy functional $\E=\E[C]$ with $\gamma=1$. By Lemma \ref{lem:convex}, $\E$ is a convex function on $\Cset$ and, therefore, the set $M$ of its minimizers,
\begin{equation}\label{eq:defM}
  M:= \left\{\tilde C\in\Cset; \; \E(\tilde C)=\min_{C\in \Cset}\E(C) \right\}
\end{equation}
is a closed convex subset of $\Cset$.
Moreover, due to the coercivity
\[
   \E[C] \geq\nu \sum_{(i,j)\in\Eset} C_{ij} L_{ij}
\]
it is nonempty.
Its structure is characterized by the following Theorem.

\begin{theorem}\label{thm:4}
Let $\gamma = 1$ and let $M$ be defined by \eqref{eq:defM}. Then the extremal points of $M$ represent loop-free graphs, and vice-versa, the loop-free elements of $M$ are extremal points of $M$.
\end{theorem}

\begin{proof}
We first prove that the extremal points of $M$ represent loop-free graphs.
Let $\tilde C\in M$. Since $\E[\tilde C] < +\infty$,
the Kirchhoff law \eqref{eq:kirchhoff} admits a solution $\tilde P\in\R^{|V|}$.
Formula \eqref{eq:derivative} gives
 \[
    \frac{\partial \E[\tilde C]}{\partial\tilde C_{ij}} = -\frac{(\tilde P_j-\tilde P_i)^2}{L_{ij}}+\nu L_{ij}.
 \]
If $\tilde C_{ij}>0$ for some $(i,j)\in E$, 
we have $\frac{\partial \E[\tilde C]}{\partial\tilde C_{ij}} =0$ and, consequently,
\begin{equation}\label{eq:thm4:1}
   {(\tilde P_j-\tilde P_i)^2} = \nu {L_{ij}^2} > 0.
\end{equation}
For contradiction, let us now assume that $\tilde C$ contains a loop, i.e., that there is a chain of edges
\begin{equation}\label{eq:thm4:2}
{\mathcal T}=\{(i_0,i_1),\dots,(i_{K-2},i_{K-1}),(i_{K-1},i_0)\}\subset E,
\end{equation}
such that $\tilde C_{i_j,i_{j+1}}>0$ for $j=0,\dots,K-1$ (we count the indices modulo $K$, so $i_{K}\equiv i_0$).
We also assume that ${\mathcal T}$ is the shortest loop, i.e., that all $i_0,\dots i_{K-1}$
are mutually different.
By \eqref{eq:thm4:1} we know that
$\tilde P_{i_j}\neq \tilde P_{i_{j+1}}$ and, with \eqref{eq:flowrate}, the corresponding fluxes $\tilde Q_{i_j,i_{j+1}}$ are all nonvanishing. 

We rewrite \eqref{eq:thm4:1} as
\begin{equation}\label{eq:thm4:3}
\tilde P_{i_j}-\tilde P_{i_{j+1}}=\sign(\tilde P_{i_j}-\tilde P_{i_{j+1}})\,\sqrt{\nu}\,L_{i_j,i_{j+1}}
\end{equation}
and sum these identities over $j=0,\dots,K-1$ to obtain
\begin{equation}\label{eq:thm4:4}
  0 = \sum_{j=0}^{K-1}\sign(\tilde P_{i_j}-\tilde P_{i_{j+1}})\,L_{i_j,i_{j+1}}.
\end{equation}
We now perturb $\tilde Q$ by adding a circular flow along ${\mathcal T}.$
For $(i,j)\in E\setminus {\mathcal T}$ we set
\[
    Q_{ij}:=\tilde Q_{ij}, \qquad
    C_{ij}:=\tilde C_{ij},
\]
while for $l=0,\dots,K-1$,
\(  \label{eq:thm4:5}
   Q_{i_l,i_{l+1}}:=\tilde Q_{i_l,i_{l+1}}+\varepsilon,\qquad
   C_{i_l,i_{l+1}}:=\tilde C_{i_l,i_{l+1}}+\frac{\varepsilon}{\sqrt{\nu}}\cdot\sign(\tilde P_{i_l}-\tilde P_{i_{l+1}}),
\)
with some $\varepsilon\neq 0$.
Obviously, 
if $|\varepsilon|$ is small enough, we have $C\in\Cset.$
Also the local mass conservation $\sum_{j\in N(i)} Q_{ij}=S_i$ is verified for all $i\in V.$ Moreover, if $(i,j)\not\in {\mathcal T}$, then $Q_{ij}=C_{ij} \frac{\tilde P_i-\tilde P_j}{L_{ij}}$ remains true.
On the other hand, for $(i_l,i_{l+1})\in{\mathcal T}$ we use \eqref{eq:thm4:3} and \eqref{eq:thm4:5} to get
\begin{align*}
   C_{i_l,i_{l+1}}\cdot \frac{\tilde P_{i_l}-\tilde P_{i_{l+1}}}{L_{i_l,i_{l+1}}}&= \left[\tilde C_{i_l,i_{l+1}}+\frac{\varepsilon}{\sqrt{\nu}}\cdot\sign(\tilde P_{i_l}-\tilde P_{i_{l+1}})\right]\cdot\frac{\tilde P_{i_l}-\tilde P_{i_{l+1}}}{L_{i_l,i_{l+1}}}\\
   &=\tilde Q_{i_l,i_{l+1}}+\frac{\varepsilon}{\sqrt{\nu}}\cdot\sign(\tilde P_{i_l}-\tilde P_{i_{l+1}})\cdot \sign(\tilde P_{i_l}-\tilde P_{i_{l+1}})\sqrt{\nu}\\
&=\tilde Q_{i_l,i_{l+1}}+\varepsilon=Q_{i_l,i_{l+1}}.
\end{align*}
Finally, using successively \eqref{eq:thm4:1}, \eqref{eq:thm4:5} and \eqref{eq:thm4:4}, we have
\begin{align*}
 \sum_{l=0}^{K-1}\left[C_{i_l,i_{l+1}}\frac{(\tilde P_{i_l}- \tilde P_{i_{l+1}})^2}{L_{i_l,i_{l+1}}^2}+\nu C_{i_l,i_{l+1}}\right]L_{i_l,i_{l+1}} 
  &= 2\nu\sum_{l=0}^{K-1}  C_{i_l,i_{l+1}}L_{i_l,i_{l+1}} \\
  &= 2\nu \sum_{l=0}^{K-1}\left[\tilde C_{i_l,i_{l+1}}+\frac{\varepsilon}{\sqrt{\nu}}\cdot\sign(\tilde P_{i_l}-\tilde P_{i_{l+1}})\right]L_{i_l,i_{l+1}}\\
 &= 2\nu\sum_{l=0}^{K-1} \tilde C_{i_l,i_{l+1}}L_{i_l,i_{l+1}}+2\varepsilon\sqrt{\nu}\sum_{l=0}^{K-1}\sign(\tilde P_{i_l}-\tilde P_{i_{l+1}})L_{i_l,i_{l+1}}  \\
 &= 2\nu\sum_{l=0}^{K-1} \tilde C_{i_l,i_{l+1}}L_{i_l,i_{l+1}},
\end{align*}
so that one more application of \eqref{eq:thm4:1} gives
\[
  \sum_{l=0}^{K-1}\left[C_{i_l,i_{l+1}}\frac{(\tilde P_{i_l}- \tilde P_{i_{l+1}})^2}{L_{i_l,i_{l+1}}^2}+\nu C_{i_l,i_{l+1}}\right]L_{i_l,i_{l+1}} 
  = 
  \sum_{l=0}^{K-1}\left[\tilde C_{i_l,i_{l+1}}\frac{(\tilde P_{i_l}-\tilde P_{i_{l+1}})^2}{L_{i_l,i_{l+1}}^2}+\nu \tilde C_{i_l,i_{l+1}}\right] L_{i_l,i_{l+1}}.
\]
This immediately implies that ${\mathcal E}[C] = {\mathcal E}[\tilde C]$.
It follows that if $\tilde C$ contains a loop, then it can be written as a convex combination of $C^{(+)}$ and, resp., $C^{(-)}\in\Cset$, which are obtained from $\tilde C$ by taking (small) positive and, resp., negative $\eps$ in \eqref{eq:thm4:5}.
Therefore, $\tilde C$ is not an extremal point of $M$. This concludes the proof that all extremal points of $M$ are loop-free.

Let us now prove the converse claim, namely, that every loop-free element of $M$ is an extremal point.
Let $C\in M$ be loop-free and assume, for contradiction, that $C$ can be written as $C=\alpha C^{(1)}+(1-\alpha)C^{(2)}$ for some $0<\alpha<1$ and
$C^{(1)},C^{(2)}\in M$. But then
\[
   \bigl\{ (i,j)\in E:C_{ij}>0 \bigr\} =
      \bigl\{(i,j)\in E:C^{(1)}_{ij}>0 \bigr\}\cup
      \bigl\{(i,j)\in E:C^{(2)}_{ij}>0 \bigr\},
\]
which means that $C^{(1)}$ and $C^{(2)}$ are also loop-free. As $\E(C)=\E(C^{(1)})=\E(C^{(2)})<+\infty$, the Kirchhoff law \eqref{eq:kirchhoff} is solvable with $C$, $C^{(1)}$ and, resp., $C^{(2)}$, and \eqref{eq:flowrate} gives the corresponding fluxes $Q$, $Q^{(1)}$ and, resp., $Q^{(2)}$. 
As for loop-free networks the fluxes are uniquely determined by the sources/sinks $S$, we have
$Q=Q^{(1)}=Q^{(2)}$ and, in turn, $C=C^{(1)}=C^{(2)}$. Hence, $C$ is an extremal point.
\end{proof}

\section{Introducing robustness}\label{sec:robustness}

The Cheeger constant (also called isoperimetric number) of a graph is a numerical measure of whether or not a graph has a ``bottleneck".
For $G=(V,E)$ its Cheeger constant $\chee(G)$ is defined as
\[
   \chee(G) = \min \left\{ \frac{|\partial W|}{|W|}; \; \emptyset \neq W \subset V, |W| \leq \frac{|V|}{2} \right\},
\]
where $\partial W \subset E$ denotes the set of edges having one end in $W$ and the other end in $V\setminus W$.

The Cheeger constant is strictly positive if and only if $G$ is a connected graph. Intuitively, if the Cheeger constant is small but positive, then there exists a ``bottleneck", in the sense that there are two ``large" sets of vertices with ``few" links (edges) between them. The Cheeger constant is large if any possible division of the vertex set into two subsets has ``many" links between those two subsets. In other words, it is a measure of 
resilience of the network against disturbance of connectivity through removal of edges.

In general, the calculation of the isoperimetric number of graphs with multiple edges is NP hard \cite{Mohar}.
Consequently, we consider a related quantity, the algebraic connectivity, also called Fiedler number,
defined as the second smallest (i.e., for a connected graph, the smallest nonzero) eigenvalue of the matrix Laplacian $\Lap = \Lap[A]$ of the adjacency matrix $A$.
The Fiedler number $\f(G)$ is another classical measure of connectivity of the graph, or robustness with respect to edge removal.
It is related to the isoperimetric number by the well known Cheeger inequalities \cite{Mohar},
\[
   \frac{\chee(G)^2}{2\Delta(G)} \leq \f(G) \leq  2 \chee(G),
\]
where $\Delta(G)$ is the maximum degree for the nodes in $G$.

In our paper we shall work with the generalized Fiedler number $\f[C]$, calculated as the second smallest eigenvalue of the weighted matrix Laplacian $\Lap[C]$ with weights given by the edge conductivities $C\in\Cset$.
This reflects the modeling assumption that edges with higher conductivities contribute more to overall network robustness (for instance, because they may be more resilient against severing).
Our idea is to modify the energy functional \eqref{eq:energy} to take into account robustness in terms of the generalized Fiedler number $\f[C]$.
Consequently, for $\mu>0$ we introduce the modified energy functional $\F=\F[C]$,
\(  \label{eq:F}
   \F[C] := \E[C] - \mu\,\minL \, \frac{|V|-1}{2} \, \f[C],
\)
where $\E=\E[C]$ is defined in \eqref{eq:energy} and
\(  \label{def:minL}
    \minL := \min_{(i,j)\in\Eset} L_{ij} >0.
\)
The reason for scaling the second term in \eqref{eq:F} by $\minL$ is that the energy \eqref{eq:energy} is homogeneous with respect to the edge lengths,
i.e., a multiplication of $L=(L_{ij})_{(i,j)\in\Eset}$ by a positive factor leads to a rescaling of the value of the energy by the same factor.
Consequently, the modified energy \eqref{eq:F} has the same scaling property.
The motivation for introducing the factor $\frac{|V|-1}{2}$ is clarified in Lemma \ref{lem:boundbelow} below.

Recalling that the model \eqref{eq:flowrate}--\eqref{eq:energy} bears relevance in biological applications as long as $\gamma\leq 1$, we consider the modified energy functional \eqref{eq:F} exclusively with $\gamma=1$ in the sequel.
The reason for excluding the values $\gamma<1$ is that $\F=\F[C]$ with $\mu>0$ is then, in general, not bounded from below:

\begin{example}\label{ex:nbb}
Let us consider a complete graph with all edges of the same length, i.e., $L_{ij} = \minL$ for all $(i,j)\in\Eset$.
Moroever, for some $c>0$, let us take $C_{ij} = c$ for all $(i,j)\in\Eset$.
Then we have $\f[C] = c |V|$ and
\[
   \F[C] &=& \Ekin[C] + \nu \sum_{(i,j)\in\Eset} C_{ij}^\gamma L_{ij} - \mu\,\minL\, \frac{|V|-1}{2}\, \f[C]  \\
      &=& \Ekin[C]  +  \minL\, |V|\, \frac{|V|-1}{2} \left(\nu c^\gamma - \mu c \right).
\]
Lemma \ref{lem:ekinbound} gives
\[
   \Ekin[C] \leq \frac{\Norm{S}^2}{\f[\widetilde C]} = \frac{\Norm{S}^2 \minL}{c |V|},
\]
with $\widetilde C_{ij} := C_{ij}/L_{ij} = c/\minL$. 
Consequently, we readily have $\lim_{c\to+\infty} \F[C] = -\infty$ whenever $\gamma<1$ and $\mu>0$.
\end{example}

On the other hand, for $\gamma=1$ we have 
the following Lemma establishing boundedness of $\F=\F[C]$ from below if $\mu\leq\nu$ and coercivity if $\mu>\nu$.
It is a direct consequence of \cite[Claim 3.5]{Fiedler}, which gives an upper bound on the value of the Fiedler number.
We provide complete proof here for the sake of the reader.

\begin{lemma}\label{lem:boundbelow}
Let $\gamma=1$ and $\mu \leq \nu$.
Then the modified energy functional $\F=\F[C]$, defined in \eqref{eq:F}, satisfies
\(  \label{eq:boundbelow}
    \F[C] \geq \Ekin[C] \geq 0 \qquad\mbox{for all } C\in\Cset.
\)
Moreover, if $\mu < \nu$, then $\F=\F[C]$ is coercive in the sense that there exists $\alpha >0$ such that
\(   \label{coercivity}
   \F[C] \geq \alpha \sum_{i=1}^{|V|} \sum_{j=1}^{|V|} C_{ij} \qquad\mbox{for all } C\in\Cset.
\)
\end{lemma}

\begin{proof}
Let $\Lap$ be the Laplacian matrix corresponding to $C\in\Cset$ and let $\f[C]$ be the Fiedler number of $C$. We use the Courant-Fisher representation \eqref{Courant} of $\f[C]$ in the form
\(  \label{Courant2}
   \f[C] = \min \left\{v^T \Lap v;\, v\in\R^{|V|}, \, \Norm{v} = 1,\, v^T\mathbf{1} = 0 \right\}.
\)
We claim that the matrix
\[
   \widetilde \Lap := \Lap - \left( I - \frac{\mathbf{1}\otimes\mathbf{1}}{|V|} \right) \f[C]
\]
is positive semidefinite.
Let us choose any vector $y\in\R^{|V|}$, then we have the decomposition $y = \alpha \mathbf{1} + \beta v$,
where $v\perp\mathbf{1}$ and $\Norm{v}=1$.
Since, trivially, $\widetilde \Lap \mathbf{1} = 0$, we have
\[
   y^T \widetilde \Lap y = \beta^2 v^T \widetilde \Lap v = \beta^2 \left(v^T \Lap v - \f[C] \right) \geq 0,
\]
where the nonnegativity follows from \eqref{Courant2}.
Indeed, $\widetilde \Lap$ is positive semidefinite, and, consequently, all its diagonal elements are nonnegative.
In particular,
\[
   \min_{i\in\Vset} \widetilde \Lap_{ii} =  \min_{i\in\Vset} \Lap_{ii} - \left(1 - |V|^{-1} \right) \f[C] \geq 0,
\]
so that
\[
     \f[C] \leq \frac{|V|}{|V| - 1} \min_{i\in\Vset} \Lap_{ii} = \frac{|V|}{|V| - 1} \min_{i\in\Vset} \sum_{j=1}^{|V|} C_{ij} \leq  \frac{1}{|V| - 1} \sum_{i=1}^{|V|} \sum_{j=1}^{|V|} C_{ij}.
\]
Therefore, using the convention $L_{ij}=C_{ij}=0$ if $(i,j)\notin \Eset$,
we finally conclude that
\[
   \F[C] &=& \Ekin[C] + \frac{\nu}{2} \sum_{i=1}^{|V|} \sum_{j=1}^{|V|} C_{ij} L_{ij} - \mu\,\minL\, \frac{|V|-1}{2}\, \f[C] \\
      &\geq& \Ekin[C] + \left(\frac{\nu}{2} \min_{(i,j)\in\Eset} L_{ij} - \frac{\mu \,\minL}{2} \right) \sum_{i=1}^{|V|} \sum_{j=1}^{|V|} C_{ij} \\
        &=& \Ekin[C]  + (\nu-\mu) \frac{\minL}{2}  \sum_{i=1}^{|V|} \sum_{j=1}^{|V|} C_{ij},
\]
with $\minL = \min_{(i,j)\in\Eset} L_{ij}>0$ as defined in \eqref{def:minL}.
For $\mu\leq\nu$ this immediately gives \eqref{eq:boundbelow}.
Moreover, if $\mu < \nu$, then \eqref{coercivity} holds with
$\alpha := (\nu-\mu) \frac{\minL}{2} >0$.
\end{proof}

\begin{remark}\label{rem:optimal}
The statement of Lemma \ref{lem:boundbelow} is optimal in the following sense:
Considering again the setting from Example \ref{ex:nbb}, we have
\[
   \F[C] &=& \Ekin[C] + \nu \sum_{(i,j)\in\Eset} C_{ij} L_{ij} - \mu\,\minL\, \frac{|V|-1}{2}\, \f[C]  \\
      &\leq& \frac{\Norm{S}^2 \minL}{c |V|}  +  \minL\, |V|\, \frac{|V|-1}{2} \left(\nu - \mu \right) c.
\]
Consequently, $\lim_{c\to+\infty} \F[C] = 0$ if $\mu = \nu$,
and $\lim_{c\to+\infty} \F[C] = -\infty$ if $\mu > \nu$.
\end{remark}

Let us note that the idea of promoting robustness of the transportation network by means of the Fiedler number $\f[C]$ can also be realized by minimizing \eqref{eq:energy} subject to \eqref{eq:flowrate}--\eqref{eq:kirchhoff} on the subset of $\Cset$ where $\f[C]\geq \mu$, for a given $\mu>0$. With this setting the values of $\gamma<1$ are perfectly admissible. We plan to explore this direction in a future work.

A fundamental property of the functional \eqref{eq:F} with $\gamma=1$ is its convexity.

\begin{lemma} \label{lem:Fconvex}
The functional \eqref{eq:F} with $\gamma=1$ is convex on the set $\Cset$.
\end{lemma}
\begin{proof}
We recast \eqref{Courant2} as
\[
   \f[C] 
   = \min \left\{\sum_{(i,j)\in\Eset} C_{ij}(v_i - v_j)^2;\, v\in\R^{|V|}, \, \Norm{v} = 1,\,
   v^T\mathbf{1} = 0\right\}.
\]
Therefore, $\f[C]$ is a minimum of linear functionals in $C$ and thus concave. The claim follows directly from the convexity of $\E=\E[C]$ with $\gamma=1$ established in Lemma \ref{lem:convex}.
\end{proof}

Let us now present two toy models, where in a very simple setting we are able to provide explicit results for the minimizers of \eqref{eq:F} with $\gamma=1$ and their dependence on the value of $\mu>0$. 

\subsection{Toy model: A triangle with sources/sinks $+1$, $-1$, $0$}\label{subsec:toy1}
Let us consider a network consisting of three nodes $\Vset = \{0,1,2\}$
and three (unoriented) edges $\Eset=\{(0,1),(0,2),(1,2)\}$ of unit length.
The sources and sinks be given by
\[
   S_0=1,\qquad S_1=-1,\qquad S_2=0.
\]
Taking into account the symmetry of the problem,
we denote $C_0$ the conductivity of the edge $(0,1)$,
and $C_1$ the conductivity of the other two edges, $(0,2)$ and $(1,2)$.
Moreover, let us denote $Q_0$ the flux through the edge $(0,1)$;
obviously, the flux through the edges $(0,2)$ and $(1,2)$ is $1-Q_0$.
The Kirchhoff law yields, by a simple calculation,
\[
   Q_0 = \frac{2C_0}{2C_0 + C_1}, \qquad 1-Q_0 = \frac{C_1}{2C_0 + C_1}.
\]
Then the energy functional \eqref{eq:energy} with $\gamma=1$ reads
\(   \label{En110}
   \mathcal{E}[C] &=& \frac{Q_0^2}{C_0} + 2\frac{(1-Q_0)^2}{C_1} + \nu\left( C_0 + 2 C_1 \right) \\
    &=& \frac{2}{2C_0 + C_1} + \nu\left( C_0 + 2 C_1 \right),
\)
with metabolic coefficient $\nu>0$.
The global minimizer $\bar C$ of $\mathcal{E}[C]$ is $\bar C_0 = 1/\sqrt\nu$, $\bar C_1=0$, with $\mathcal{E}[\bar C] = 2\sqrt\nu$.
Obviously, with $\bar C_1=0$, only the edge $(0,1)$ is present and, consequently, the graph is disconnected.
However, since $S_2=0$, the Kirchhoff law \eqref{eq:kirchhoff} is solvable and network still fulfills its task
of transporting unit mass from node $0$ to node $1$.

\begin{figure} \centering
 \includegraphics[width=.47\textwidth]{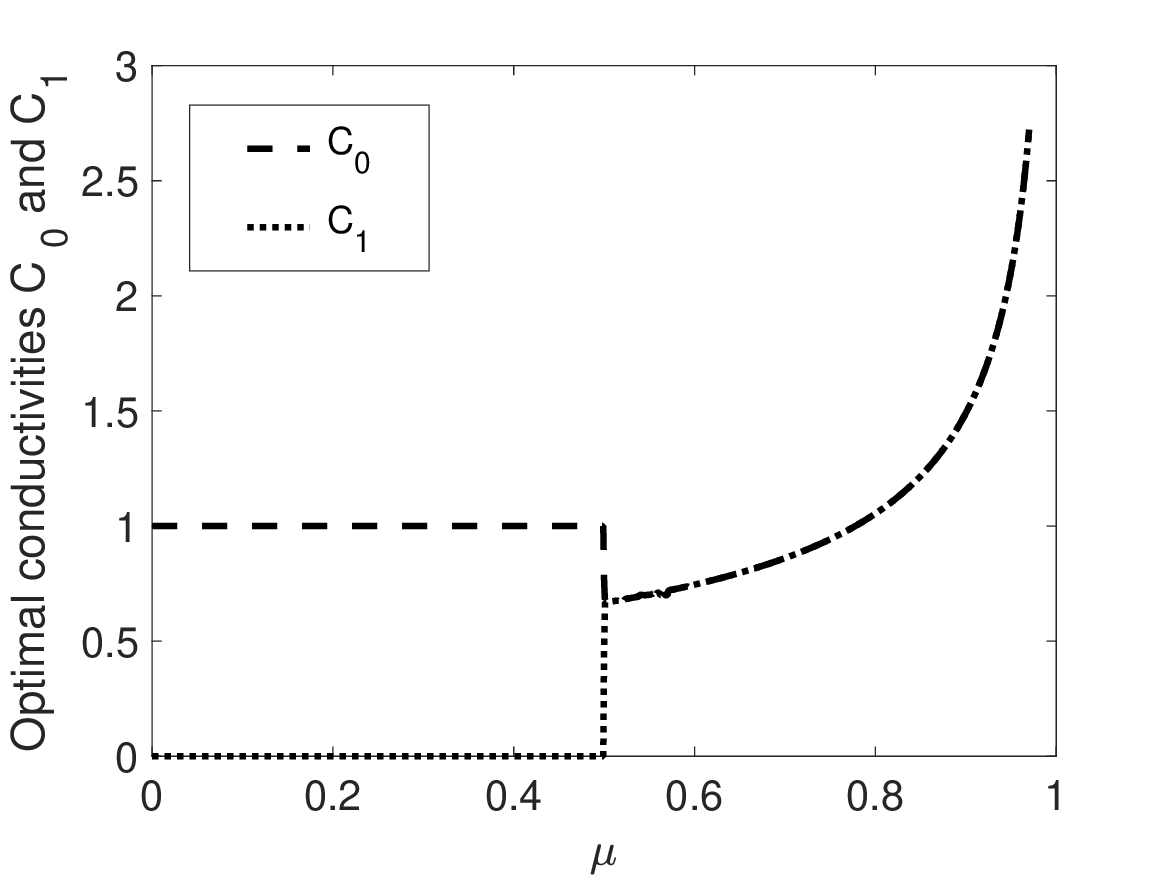}
 \includegraphics[width=.47\textwidth]{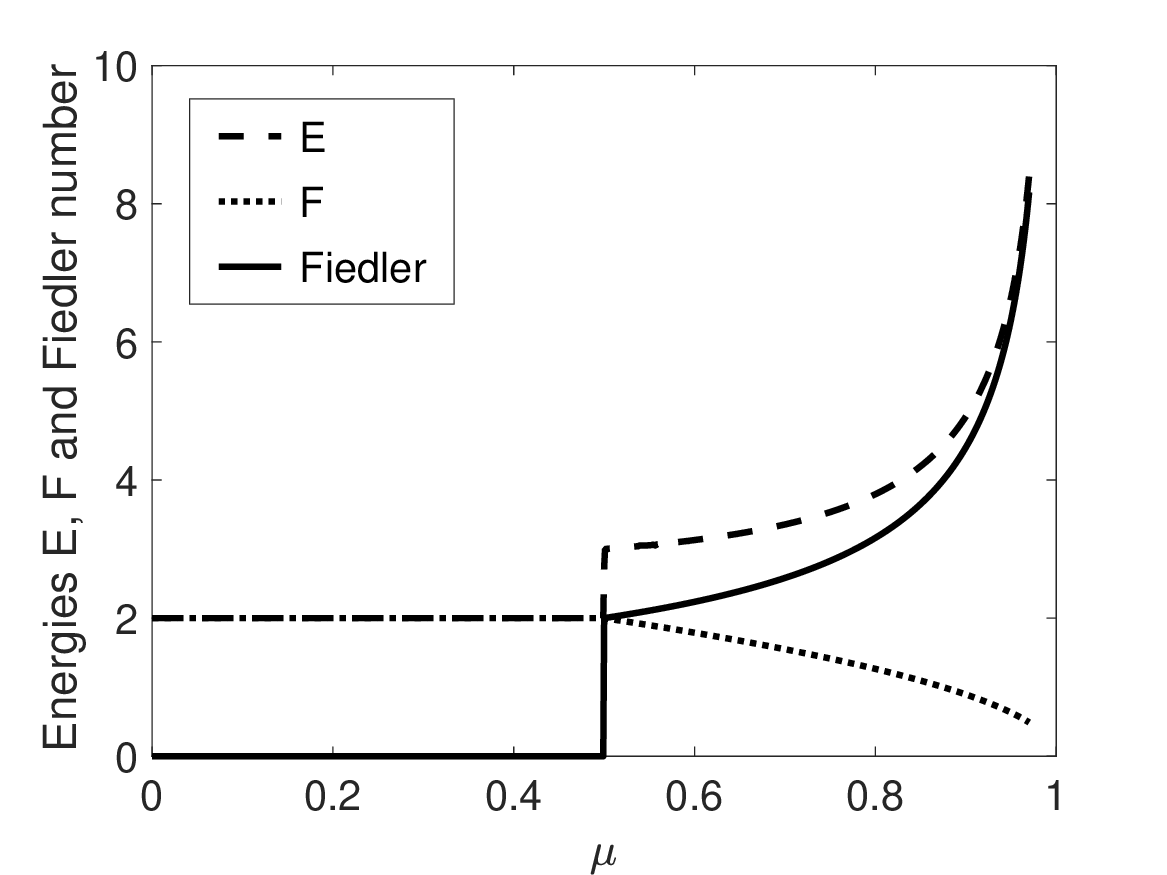}
 \caption{Left panel: Optimal values of the conductivities $C_0$ and $C_1$ for the minimization problem \eqref{Fn110} with $\nu:=1$, as a function of the parameter $\mu\in [0,1]$.
 Right panel: Values of the functionals $\E=\E[C]$, $\F=\F[C]$ and the Fiedler number $\f=\f[C]$ for the optimal solutions.
 \label{fig:T110}}
\end{figure}

The eigenvalues of the matrix Laplacian of the corresponding weighted adjacency matrix are
\[
   \lambda_0 = 0, \qquad \lambda_1 = 3C_1, \qquad \lambda_2 = 2C_0 + C_1.
\]
Consequently, the Fiedler number is $\f[C] = \min\{ 2C_0 + C_1, 3C_1 \}$ and the modified energy \eqref{eq:F} reads
\(   \label{Fn110}
   \F[C] = \frac{2}{2C_0 + C_1} + \nu\left( C_0 + 2 C_1 \right) - \mu \min\{ 2C_0 + C_1, 3C_1 \},
\)
with $\mu \geq 0$.
We note that due to the concavity of $\f[C] = \min\{ 2C_0 + C_1, 3C_1 \}$ and strict (but non-uniform) convexity
of the kinetic energy term $\frac{2}{2C_0 + C_1}$, the functional $\F[C]$ is non-uniformly strictly convex on the positive quadrant $\R^2_+$.
Lemma \ref{lem:boundbelow} states that $\F[C]$ is bounded from below if $\mu\leq\nu$ and coercive if $\mu<\nu$.
By inspecting the case $C_0=C_1$, we see that this is also a necessary condition.
In fact, for $\mu=\nu$, choosing $C_0=C_1=c$ with $c>0$, we have $\F[C] = \frac{2}{3c}$
and, consequently, $\F=\F[C]$ does not have a minimizer. We thus study the minimizers of \eqref{Fn110} for $0 < \mu < \nu$.
The result of this simple exercise, where we chose $\nu:=1$, is displayed in Fig. \ref{fig:T110}.
We observe that for $\mu<1/2$, the optimal conductivities $\bar C$ are $\bar C_0 = 1/\sqrt{\nu}$, $\bar C_1=0$,
i.e., the same as for $\mu=0$. Therefore, choosing $\mu<1/2$ does not enforce any improvement of the robustness
of the network. On the other hand, for $\mu >1/2$, the optimal conductivities are
\[
   \bar C_0 = \bar C_1 = \frac13 \sqrt{\frac{2}{\nu-\mu}}.
\]
In the right panel of Fig. \ref{fig:T110} we observe that the modified optimal energy $\F[\bar C]$
is continuous, constant for $0 < \mu < 1/2$ and decaying to zero as $\mu\to 1-$.
The optimal kinetic-metabolic energy $\E[\bar C]$ is also constant for $0 < \mu < 1/2$,
but has a discontinuity at $\mu=1/2$ and increases to $+\infty$ as $\mu\to 1-$.
This increasing branch of $\E[\bar C]$ can be interpreted as the energetic cost of enforcing the robustness of the network.

\subsection{Toy model: A triangle with sources/sinks $+1$, $-1/3$, $-2/3$}\label{subsec:toy2}
We consider a network consisting of three nodes $\Vset = \{0,1,2\}$
and three (unoriented) edges $\Eset = \{(0,1),(0,2),(1,2)\}$ of unit length
and conductivities $C_{01}$, $C_{02}$ and, resp., $C_{12}$.
The sources and sinks be given by
\[
   S_0=1,\qquad S_1=-1/3,\qquad S_2=-2/3.
\]
Let us denote by $Q=Q[C]$ the flux from vertex $0$ to vertex $1$.
Then the Kirchhoff law \eqref{eq:kirchhoff} gives
\(   \label{toyQ}
   Q = \frac{C_{01}(C_{02}+3C_{12})}{3(C_{01}C_{02}+C_{01}C_{12}+C_{02}C_{12})}.
\)
Moreover, due to the mass conservation in vertex $1$, the flux
from $0$ to $2$ equals to $1-Q$, and
the flux from vertex $1$ to $2$ is $Q-1/3$.

The energy functional \eqref{eq:energy} with $\gamma=1$ reads
\(  \label{En133}
    \E[C] = \frac{Q^2}{C_{01}} + \frac{(1-Q)^2}{C_{02}} + \frac{(Q-1/3)^2}{C_{12}} + \nu \left( C_{01} + C_{02} + C_{12} \right).
\)
An easy exercise reveals that the energy $\E[C]$ is globally minimized for
$\bar C_{01}=1/3$, $\bar C_{02} = 2/3$ and $\bar C_{12} = 0$, with $Q=1/3$.

\begin{figure} \centering
 \includegraphics[width=.47\textwidth]{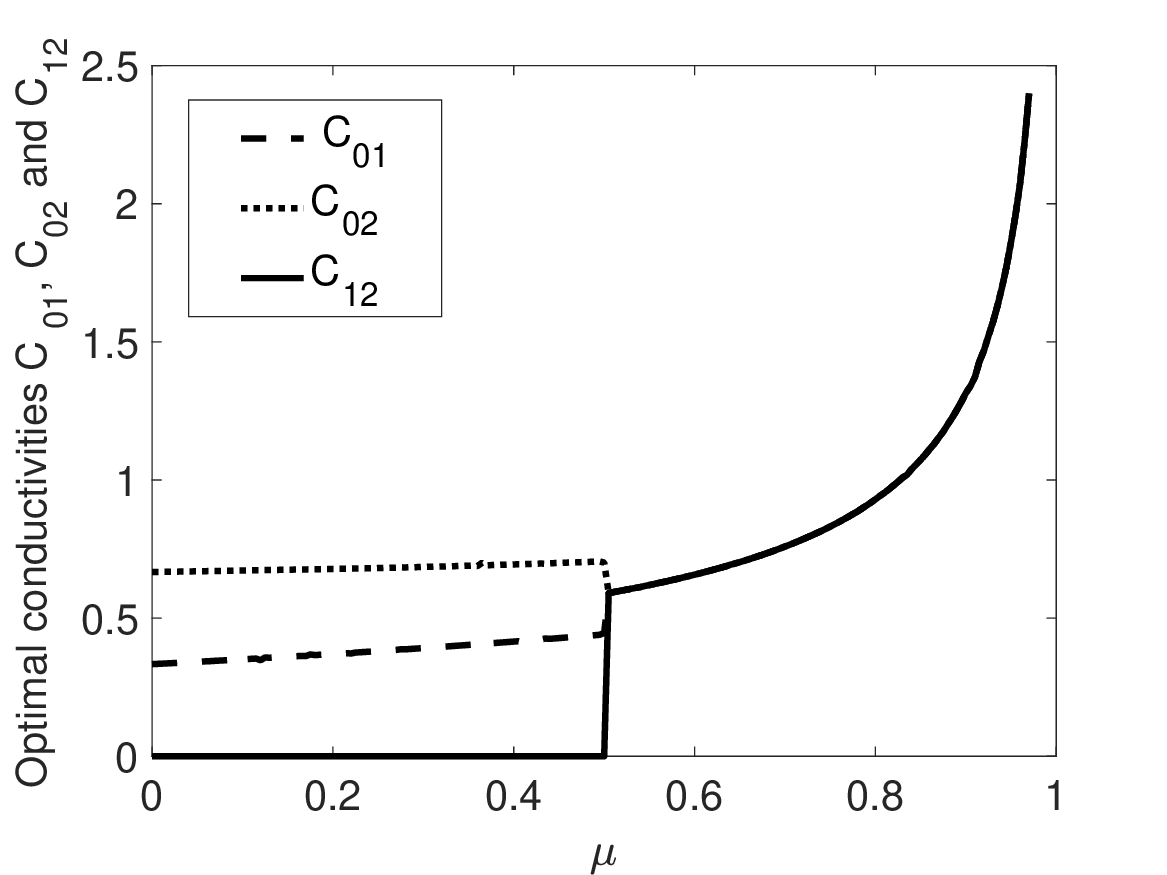}
 \includegraphics[width=.47\textwidth]{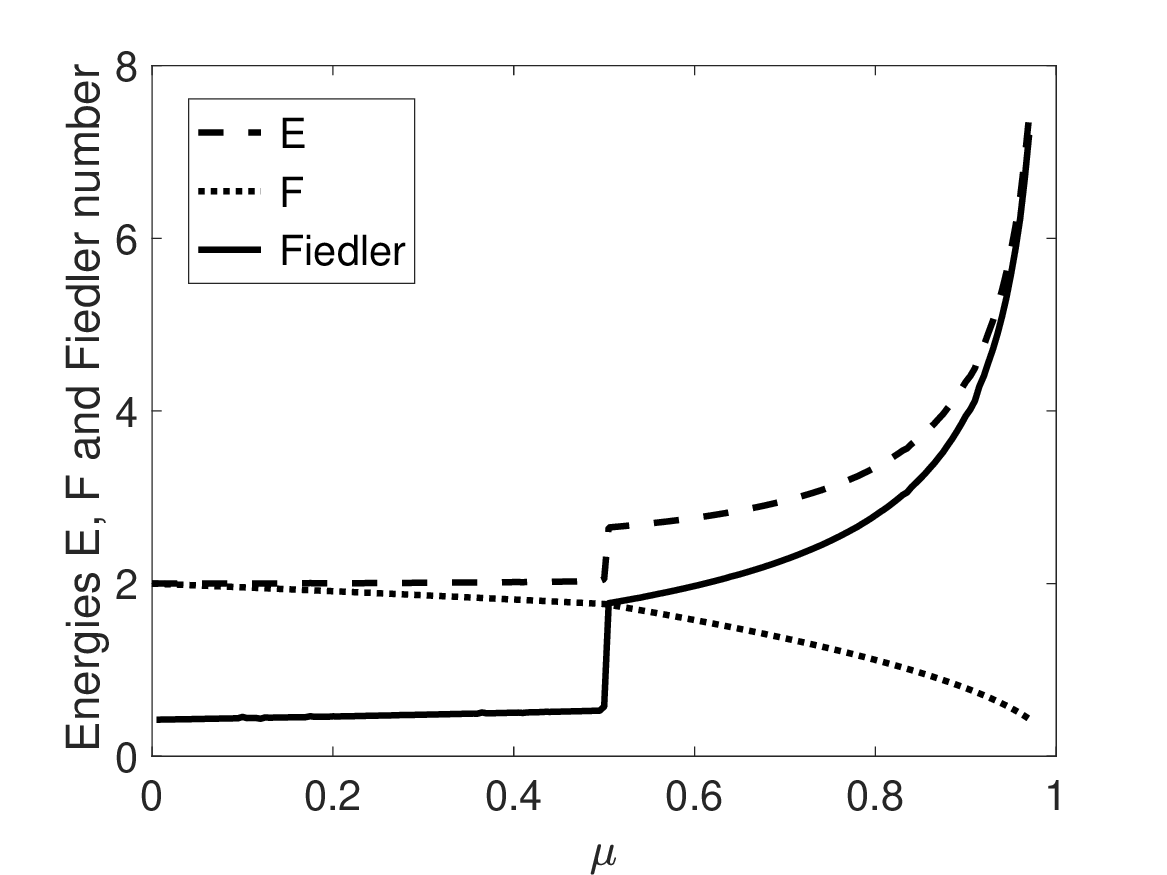}
 \caption{Left panel: Optimal values of the conductivities $C_{01}$, $C_{02}$ and $C_{12}$ for the minimization problem \eqref{eq:F133} as a function of the parameter $\mu\in [0,1]$.
 Right panel: Values of the functionals $\E=\E[C]$, $\F=\F[C]$ and the Fiedler number $\f=\f[C]$ for the optimal solutions.
 \label{fig:T133}}
\end{figure}

The eigenvalues of the matrix Laplacian of the corresponding weighted adjacency matrix
are $\lambda_0=0$ and
\[
   \lambda_1 &=& C_{01} + C_{02} + C_{12} - \sqrt{C_{01}^2 + C_{02}^2 + C_{12}^2 - C_{01}C_{02} - C_{01}C_{12} - C_{02}C_{12}},\\
   \lambda_2 &=& C_{01} + C_{02} + C_{12} + \sqrt{C_{01}^2 + C_{02}^2 + C_{12}^2 - C_{01}C_{02} - C_{01}C_{12} - C_{02}C_{12}}. 
\]
Consequently, the Fiedler number is $\f[C] = \lambda_1$.
Since $\minL=1$ and $|V|-1=2$, the functional \eqref{eq:F} takes the form
\(   \label{eq:F133}
   \F[C] = \E[C] - \mu \lambda_1,
\)
with $\E[C]$ given by \eqref{En133}.
We minimize $\F[C]$ with respect to the variables $C_{01}$, $C_{02}$, $C_{12}\geq 0$, with $Q=Q[C]$ given by \eqref{toyQ}.
The result of this optimization problem with $\nu:=1$, obtained by an application of the Matlab function {\tt fminsearch},
is displayed in Fig. \ref{fig:T133}.
In the left panel we observe that for $\mu<1/2$ we have positive optimal conductivities $\bar C_{01} \neq \bar C_{02}$,
while $\bar C_{12}=0$. In this regime the optimal $Q$ remains equal to $1/3$.
Consequently, similarly as in the example in Section \ref{subsec:toy1}, for $\mu<1/2$ the connectivity
of the network is not improved by ``activation" of the edge $(1,2)$.
However, 
the Fiedler number $\f[\bar C]$ is slowly increasing for $0\leq\mu\leq 1/2$,
from $\f[\bar C] = 1 - \sqrt{1/3} \approx 0.423$ for $\mu=0$, to $\f[\bar C] \approx 0.527$ for $\mu=1/2$. Let us also observe that $\lambda_1 < \lambda_2$, i.e., the Fiedler number is a simple eigenvalue of the matrix Laplacian.

On the other hand,
for $\mu>1/2$ we have the minimizer $\bar C_{01}=\bar C_{02}=\bar C_{12}=:c >0$, so that the modified energy functional reads
\[
   \F[c] = \frac{Q^2 + (1-Q)^2 + (Q-1/3)^2}{c} + 3 (\nu-\mu) c,
\]
and has the global minimizer $c = \frac19 \sqrt{ \frac{14}{\nu-\mu} }$, $Q=4/9$.
Consequently, for $\mu>1/2$ the network robustness is enforced by ``activation" of the edge $(1,2)$.
It is interesting to observe that, despite the ``nonsymmetry" of the values of $S_i$,
all three edges have the same optimal conductivity.
Moreover, we obviously have $\lambda_1=\lambda_2$, i.e., the Fiedler number is a double eigenvalue of the matrix Laplacian.

\section{Numerical minimization by the projected subgradient method}\label{sec:num}

In this section we present a projected subgradient method \cite{Shor} for minimization of the modified energy functional $\F=\F[C]$ given by \eqref{eq:F} with $\gamma=1$, constrained by the Kirchhoff law \eqref{eq:flowrate}--\eqref{eq:kirchhoff}.
The reason for choosing the subgradient method is that the functional is convex (by Lemma \ref{lem:convex}), and Lipschitz continuous - this follows from Lemma \ref{lemma:pert} of the Appendix,
combined with the classical Hoffman-Wielandt inequality \cite{Hoffman-Wielandt} for the Lipschitz continuity of the eigenvalues of a normal matrix.
Then, by Rademacher theorem, the functional is almost everywhere differentiable. In fact, the Fiedler number $\f=\f[C]$ is differentiable (with respect to elements $C_{ij}$ of $C$) in points $C$ where it is a simple eigenvalue of the matrix Laplacian $\Lap[C]$,
see, e.g., \cite{Rellich, Kato, Stewart-Sun}.

On the other hand,
if $\f[C]$ is a multiple eigenvalue of $\Lap[C]$, then it does not admit a derivative in classical sense. However, it is well known that the Clarke subdifferential of the function which maps a symmetric
matrix to its $m$-th smallest eigenvalue can be explicitly calculated. The subdifferential is 
identical for all choices
of the index $m$ corresponding to equal eigenvalues, and,
moreover, it coincides with the and Michel-Penot subdifferential \cite{HUL}.
Assuming that $C\in\Cset$ represents a connected graph, the Fiedler number $\f[C]$ is the second smallest eigenvalue of the matrix Laplacian $\Lap[C]$, and we are interested in calculating (any element of) its subdifferential with respect to the symmetric matrix $C$.
As the mapping $C\mapsto \Lap[C]$ is analytic, we make use of the result provided in \cite{Sun}.

\begin{lemma}\label{lem:subdiff}
Let $C\in\Cset$ represent a connected graph and let $\f[C]$ be its Fiedler number of multiplicity $r\geq 1$, i.e., the $r$-fold eigenvalue of the matrix Laplacian $\Lap[C]$.
Let $\Theta\subset\R^{|V|\times|V|}$
be the (Clarke) subdifferential of 
$\f[C]$ at $C$ and let the unit vector $v\in\R^{|V|}$ be any element of the eigenspace of $\Lap[C]$ corresponding to $\f[C]$.
Then, denoting $\mathcal{V}[v]_{ij} := (v_i-v_j)^2$ for all $i,j\in\Vset$, we have
\(   \label{eq:subdiff}
    \mathcal{V}[v] \in \Theta.
\)
\end{lemma}

\begin{proof}
We apply \cite[Theorem 3.7]{Sun} - in particular, formula (3.13), adapted to our notation, reads
\[
   \Theta = \mbox{co}
   \left\{ A\in\R^{|V|\times |V|};\,
   A_{ij} = u^T W^T \part{\Lap[C]}{C_{ij}} W u\, \mbox{ for all } i,j\in V,\, u\in\R^r, |u|=1
   \right\}.
\]
Here $W\in\R^{|V|\times r}$ denotes the matrix of column orthonormal  basis vectors of the eigenspace of 
$\Lap[C]$ corresponding to $\f[C]$. Without loss of generality we choose $v$ to be its first column.
Moreover, in the partial derivative $\part{\Lap[C]}{C_{ij}}\in\R^{|V|\times|V|}$, the symmetry of $C$ is taken into account, i.e., for $i\neq j$ it quantifies the sensitivity of $\Lap[C]$ with respect to changes in both $C_{ij}$ and $C_{ji}$.
A trivial calculation gives then, for $i\neq j$ and $k\neq m$,
\[
 \part{\Lap_{km}}{C_{ij}} = -\delta_{k,i}\delta_{m,j} 
 -\delta_{k,j}\delta_{m,i}, \qquad 
    \part{\Lap_{kk}}{C_{ij}} = \delta_{k,i} + \delta_{k,j},
\]
where we use the shorthand notation $\Lap_{km}$ for the $(k,m)$-element of $\Lap[C]$. As $\Lap[C]$ does not depend on the diagonal elements of $C$, we have $\part{\Lap[C]}{C_{ii}}=0$ for all $i\in\Vset$.
We then easily calculate, for all $i, j\in\Vset$,
\[
   A_{ij} = \sum_{s=1}^r \sum_{\sigma=1}^r \Bigl( W_{is} W_{i\sigma} + W_{js}W_{j\sigma} - W_{is}W_{j\sigma} - W_{js}W_{i\sigma} \Bigr) u_s u_\sigma.
\]
Finally, choosing $u:=(1,0,\dots,0)^T$, we get $A_{ij} = v_i^2 + v_j^2 - 2v_iv_j$ and we conclude that $\mathcal{V}[v] \in \Theta.$
\end{proof}

Obviously, Lemma \ref{lem:subdiff}
does apply also to the case when $\f[C]$ is a simple eigenvalue of $\Lap[C]$. The Fiedler number is then differentiable in classical sense and we have
\[
   \part{\f[C]}{C_{ij}} = (v_i - v_j)^2 \qquad\mbox{for all } i,j\in\Vset,
\]
where $v\in\R^{|V|}$ is the corresponding normalized eigenvector.
For interested readers we provide an elementary proof of this claim in Section \ref{subsec:derivative} of the Appendix.

We now collected all the ingredients needed for establishing the projected subgradient method for minimization of the functional \eqref{eq:F}.
We initialize the method by choosing $C^{(0)}\in\Cset$ such that the Kirchhoff law \eqref{eq:kirchhoff} with $C^{(0)}$ admits a solution.
In our practical realization, we draw the values for the elements $C^{(0)}_{ij}$, $(i,j)\in\Eset$,
randomly from the uniform distribution on the interval $(0,1)$. The elements $C^{(0)}_{ij}$ for $(i,j)\notin\Eset$ are all set to zero.
Then, we fix some $K\in\N$ and for $k=0,\dots,K$, we perform the subgradient step
\(  \label{eq:subgrad_step}
   C^{(k+1/2)}_{ij} = C_{ij}^{(k)} + \tau_k \left(\frac{\left(P^{(k)}_i - P^{(k)}_j\right)^2}{L_{ij}} - \nu L_{ij} + \mu\,\minL \, \frac{|V|-1}{2} \, \left(v^{(k)}_i - v^{(k)}_j\right)^2 \right)
\)
for all $(i,j)\in\Eset$.
Here we used the formula \eqref{eq:derivative} for the derivative of the kinetic term of the energy, and $P^{(k)}$ denotes any solution of the Kirchhoff law 
\eqref{eq:flowrate}--\eqref{eq:kirchhoff} with $C:=C^{(k)}$. Moreover, $v^{(k)}$ is any normalized eigenvector of the matrix Laplacian $\Lap[C^{(k)}]$.
We use the diminishing step size $\tau_k>0$, in particular, $\tau_k := \tau_0 / \sqrt{k}$ for some $\tau_0>0$.
The subgradient step is followed by the projection step
\( \label{eq:proj_step}
   C^{(k+1)} = \mathbb{P}_{\Cset}\left[ C^{(k+1/2)} \right].
\)
Here $\mathbb{P}_{\Cset}: \R^{|V|\times |V|}_{\mathrm{sym}} \to \Cset$ denotes the projection onto the set $\Cset$ and is realized by simply trimming the negative elements to zero,
\[
  \mathbb{P}_{\Cset}[C]_{ij} := \max\{C_{ij},0\} \qquad\mbox{for all } i,j\in\Vset.
\]
After performing the projection step, we check whether $C^{(k+1)}$ remained in the domain of $\F=\F[C]$, which is equivalent to the solvability of the Kirchhoff law \eqref{eq:flowrate}--\eqref{eq:kirchhoff} with $C:=C^{(k+1)}$. Obviously, due to continuity, $\F[C^{(k)}] < +\infty$ implies the same for $C^{(k+1)}$
for small enough step size $\tau_k>0$.
In practical numerical realization of the method, we interrupt the calculation if we detect that the linear system \eqref{eq:flowrate}--\eqref{eq:kirchhoff} became badly conditioned, and restart with a reduced $\tau_0>0$.

The subgradient method is not a descent method, i.e., it is not guaranteed that $\F[C^{(k+1)}] \leq \F[C^{(k)}]$. Therefore, we set
\[
   \widetilde C^{(K)} := \mbox{argmin}_{k=0,\dots,K} \F[C^{(k)}].
\]
Then, due to the convexity and Lipschitz continuity of the functional $\F=\F[C]$ on its domain, the standard theory \cite{Shor} provides convergence of the sequence $(\widetilde C^{(K)})_{K>0}$ to a global minimizer. Its existence follows from the coercivity of $\F=\F[C]$, Lemma \ref{lem:boundbelow}, as long as $\mu<\nu$.

We apply the projected subgradient method \eqref{eq:subgrad_step}--\eqref{eq:proj_step} to search for optimal transportation networks in two cases - a complete graph with $7$ nodes, and a leaf-shaped graph with $122$ nodes and $323$ edges. We study how the resulting optimal networks depend on the value of the parameter $\mu\geq 0$. In particular, we are interested in the number of active edges (i.e., number of nonzero elements $C_{ij}$), the multiplicity of the Fiedler number and the values of the kinetic, metabolic and modified energies.

\subsection{Example with $|V|=7$}\label{subsec:N7}
We applied the projected subgradient method \eqref{eq:subgrad_step}--\eqref{eq:proj_step}
for minimization of the functional \eqref{eq:F} with $\gamma:=1$ and $\nu:=1$, on a complete graph $G=(\Vset, \Eset)$ with $7$ nodes located in $\R^2$,
\[
   \Vset = \{1,2,\ldots, 7\}, \qquad \Eset= \{ (i,j) = (j,i);\; i, j\in \Vset \}.
\]
The node locations and the
 source/sink intensities are specified in Table \ref{Tab:1}.
 
\begin{table}[h]  \centering
\begin{tabular}{ c | c | c | c}
  node & $x$ & $y$ & $S$ \\
  \hline
  1 & 0.100 & 0.000 & 0.164 \\
  2 & 0.874 & -0.104 & 0.794 \\
  3 & 0.581 & 0.790 & -0.128 \\
  4 & -0.043 & 1.0547 & 0.936 \\
  5 & -0.945 & 0.371 & -0.299 \\
  6 & -0.818 & -0.161 & 0.750 \\
  7 & -0.206 & -1.023 & -2.217
\end{tabular}
  \caption{Locations $(x,y)\in\R^2$ of graph nodes and source/sink intensities.}
  \label{Tab:1}
\end{table}

We chose $\tau_0:=10^{-1}$ and $K:=10^6$.
To check for consistency of the  method, we ran the calculation multiple times for every fixed value of $\mu\geq 0$, each time with a different initial $C^{(0)}\in\Cset$, with elements drawn from the uniform random distribution on $(0,1)$. In all these runs (with a fixed value of $\mu\geq 0$) the method converged to the same minimizer, up to a relative error of the order $10^{-12}$.

\begin{figure}[h!] \centering
 \includegraphics[width=.32\textwidth]{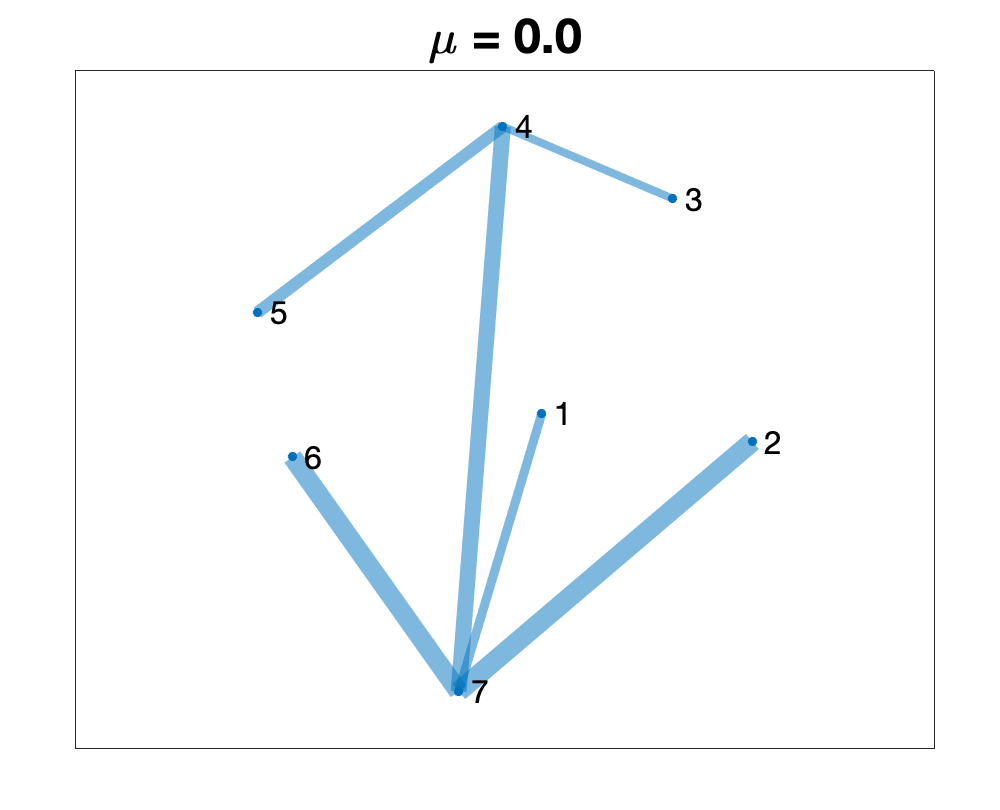}
 \includegraphics[width=.32\textwidth]{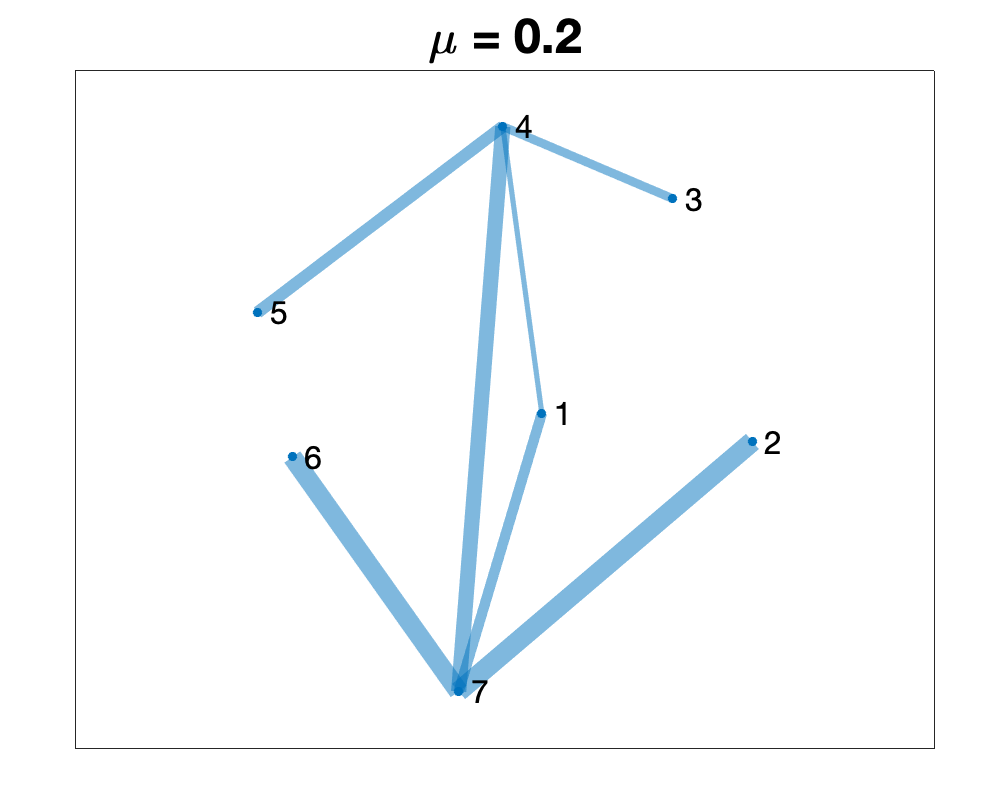}
 \includegraphics[width=.32\textwidth]{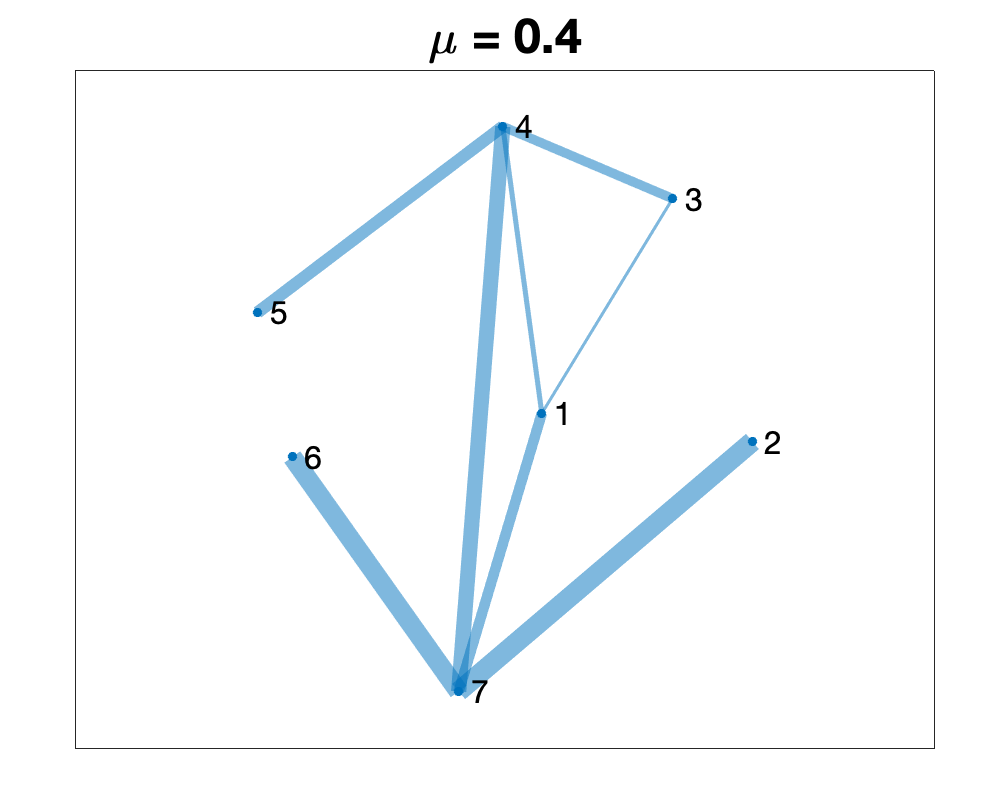} \\
 \includegraphics[width=.32\textwidth]{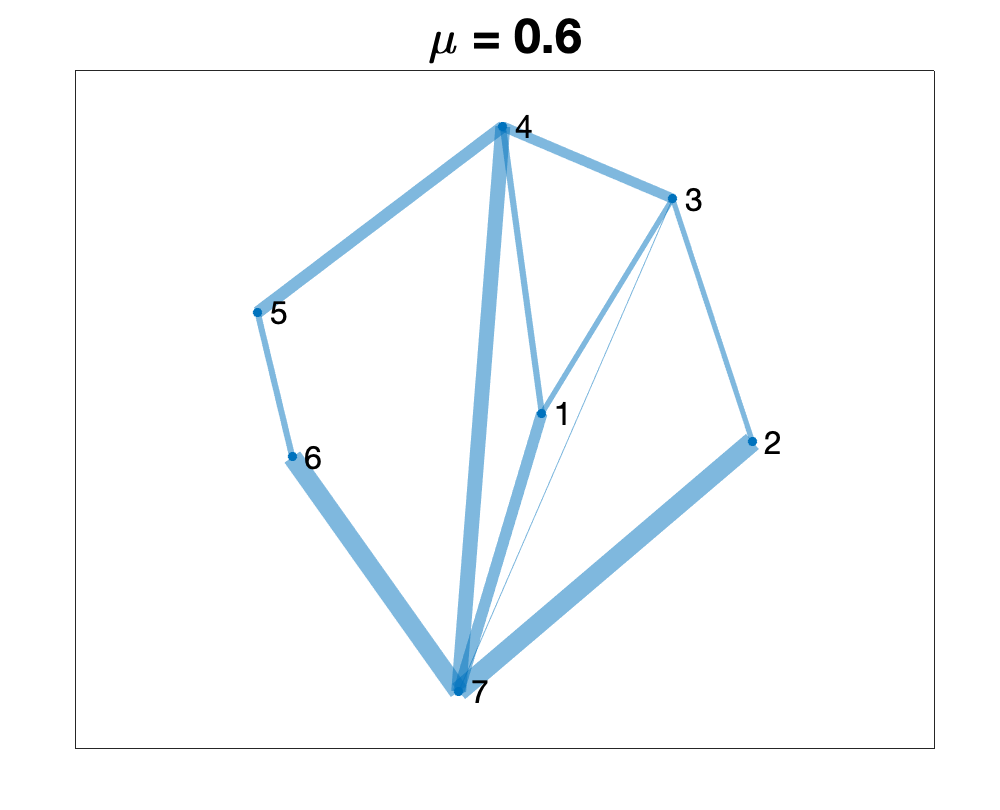}
 \includegraphics[width=.32\textwidth]{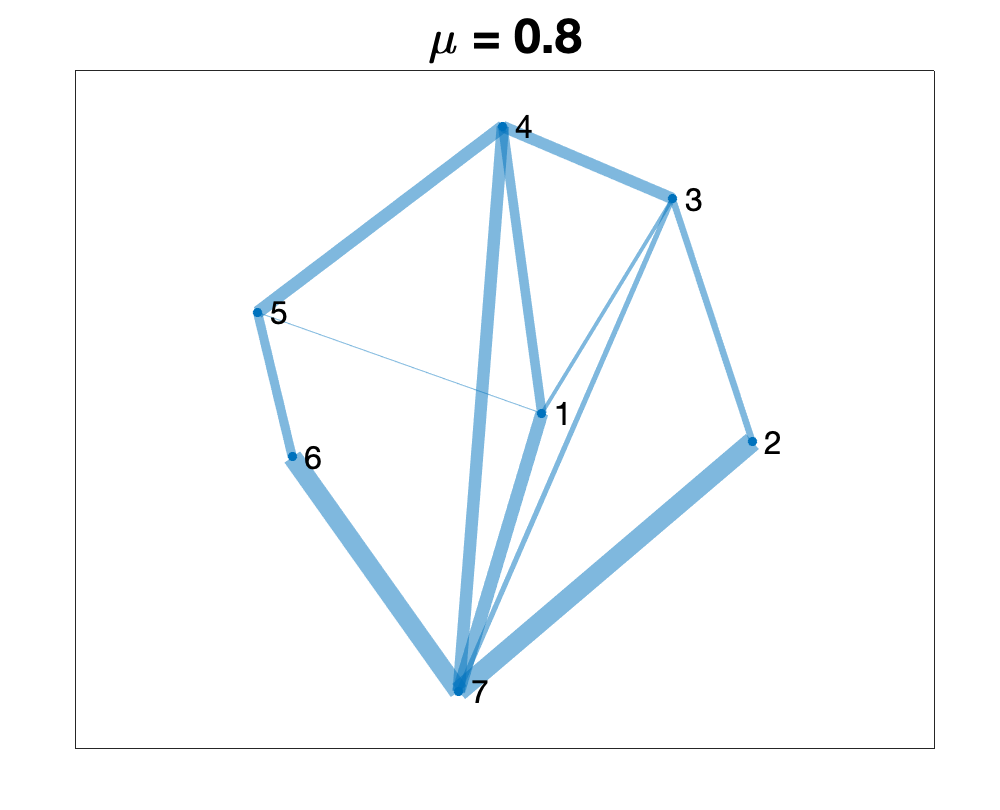}
 \includegraphics[width=.32\textwidth]{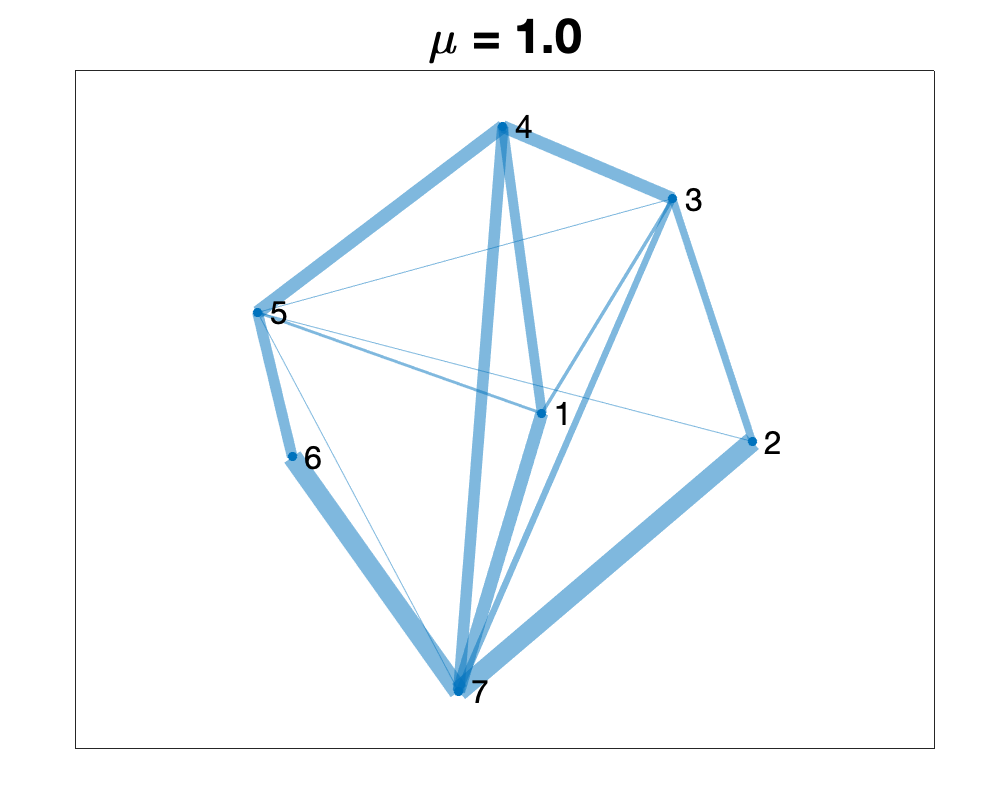} \\
 \caption{Results of minimization of the functional $\F=\F[C]$, given by \eqref{eq:F}, for the graph with $7$ vertices (Table \ref{Tab:1}) and $\mu\in\{0, 0.2, 0.4, 0.6, 0.8, 1.0\}$.
 The thickness of the line segments is proportional to the square root of the conductivity $C_{ij}$ of the corresponding edge. Edges with $C_{ij}=0$ are excluded from the plot.
 \label{fig:N7}}
\end{figure}

 The graphs corresponding to the minimizers for $\mu\in\{0, 0.2, 0.4, 0.6, 0.8, 1.0\}$ are plotted in Fig. \ref{fig:N7},
where the thickness of the line segments is proportional to the square root of the conductivity $C_{ij}\geq 0$ of the corresponding edge. Edges with $C_{ij}=0$
are excluded from the plot. We observe that for $\mu=0$ the optimal transportation structure is loop-free, which corresponds
to the result of Theorem \ref{thm:4}.
For $\mu=0.2$ one loop is present, consisting of nodes $\{1,4,7\}$.
With increasing value of $\mu$, the graph is successively becoming denser.

\begin{figure}[h!] \centering
 \includegraphics[width=.35\textwidth]{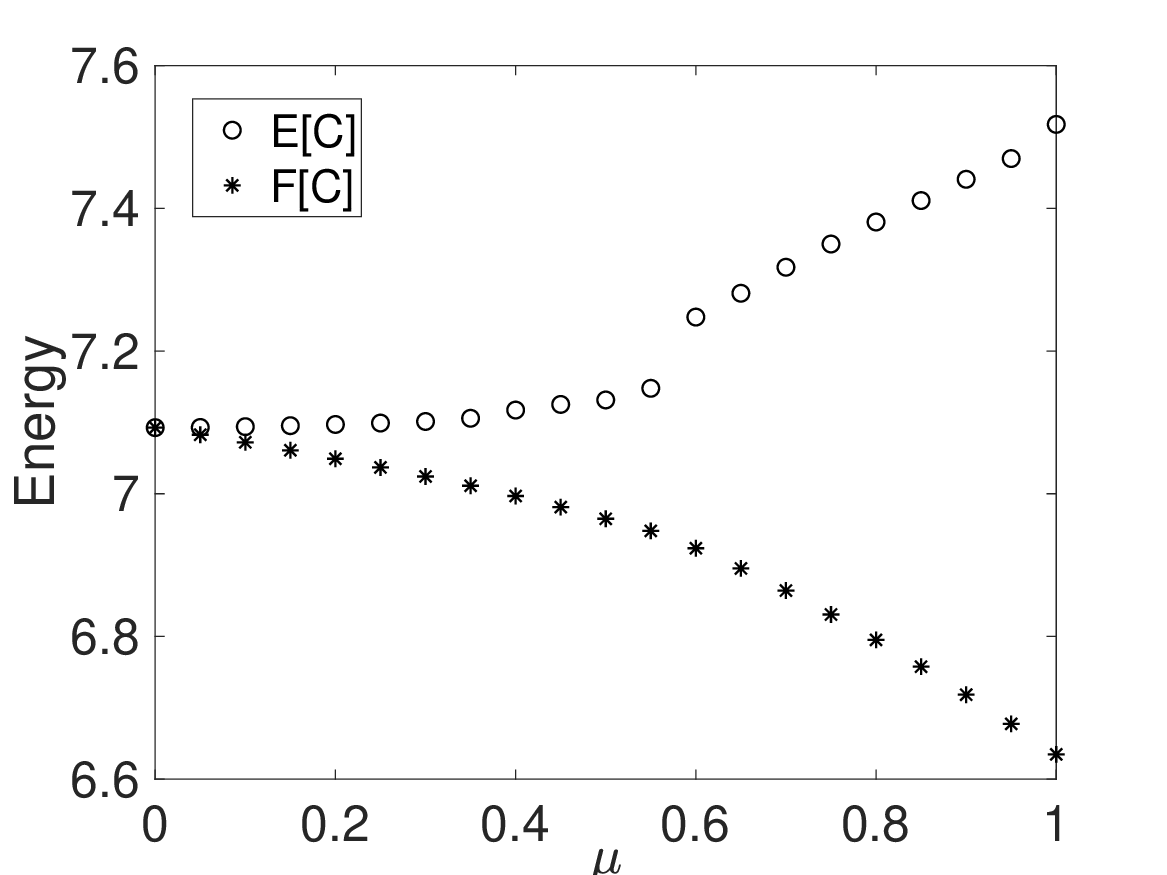} $\qquad$
 \includegraphics[width=.35\textwidth]{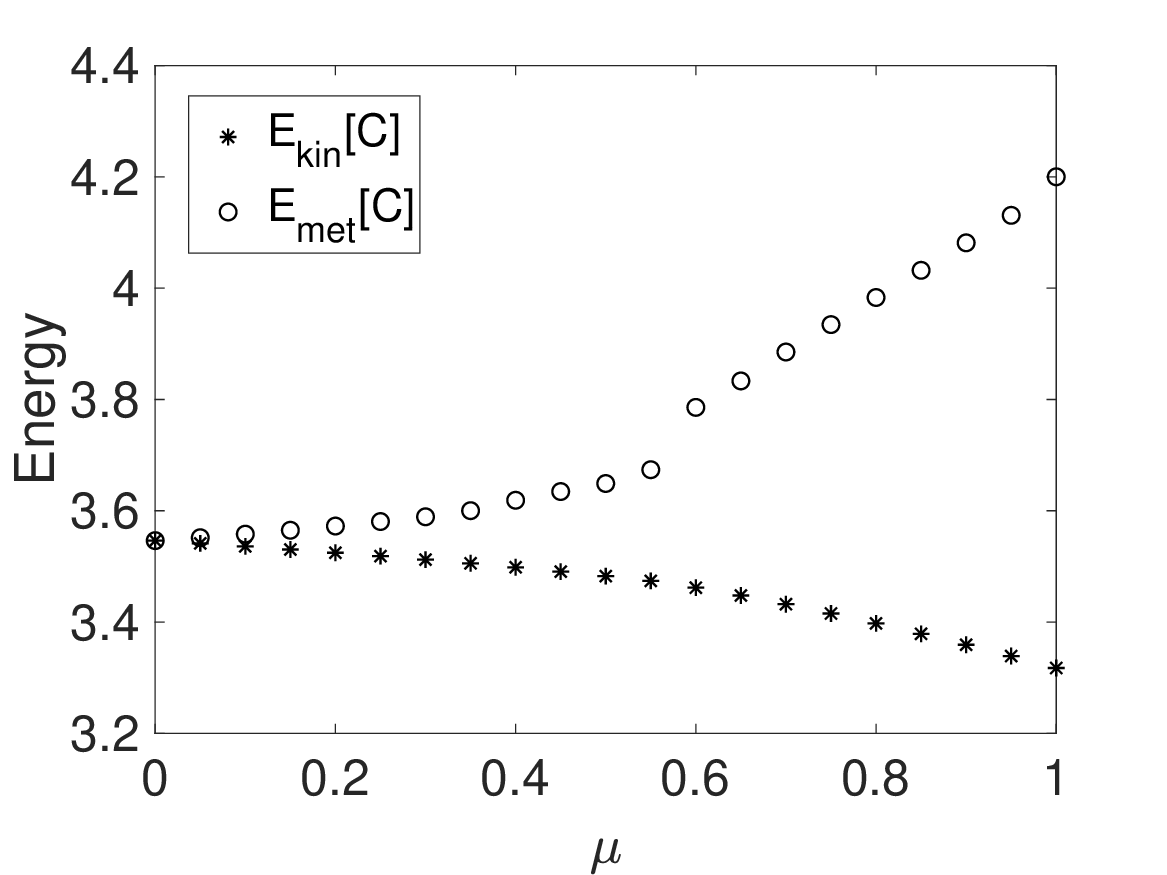} \\
 \includegraphics[width=.35\textwidth]{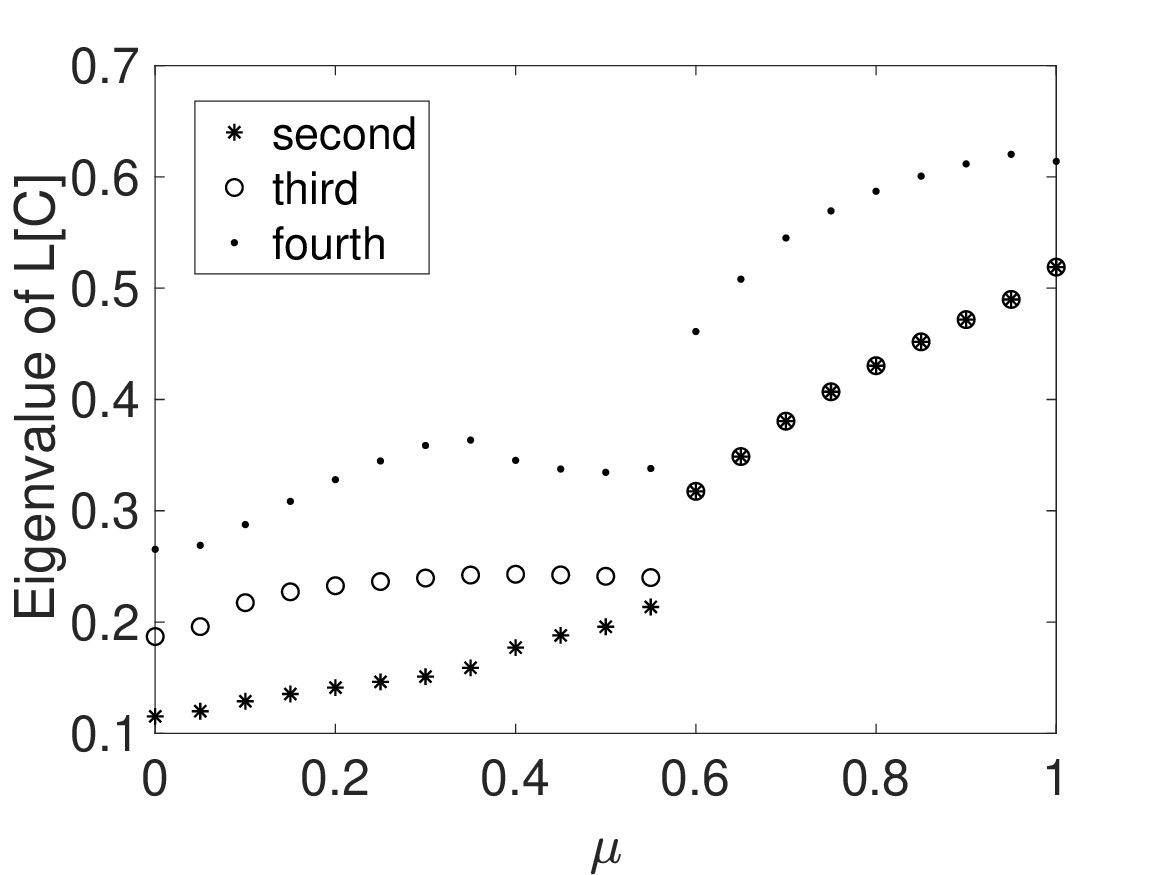} $\qquad$
 \includegraphics[width=.35\textwidth]{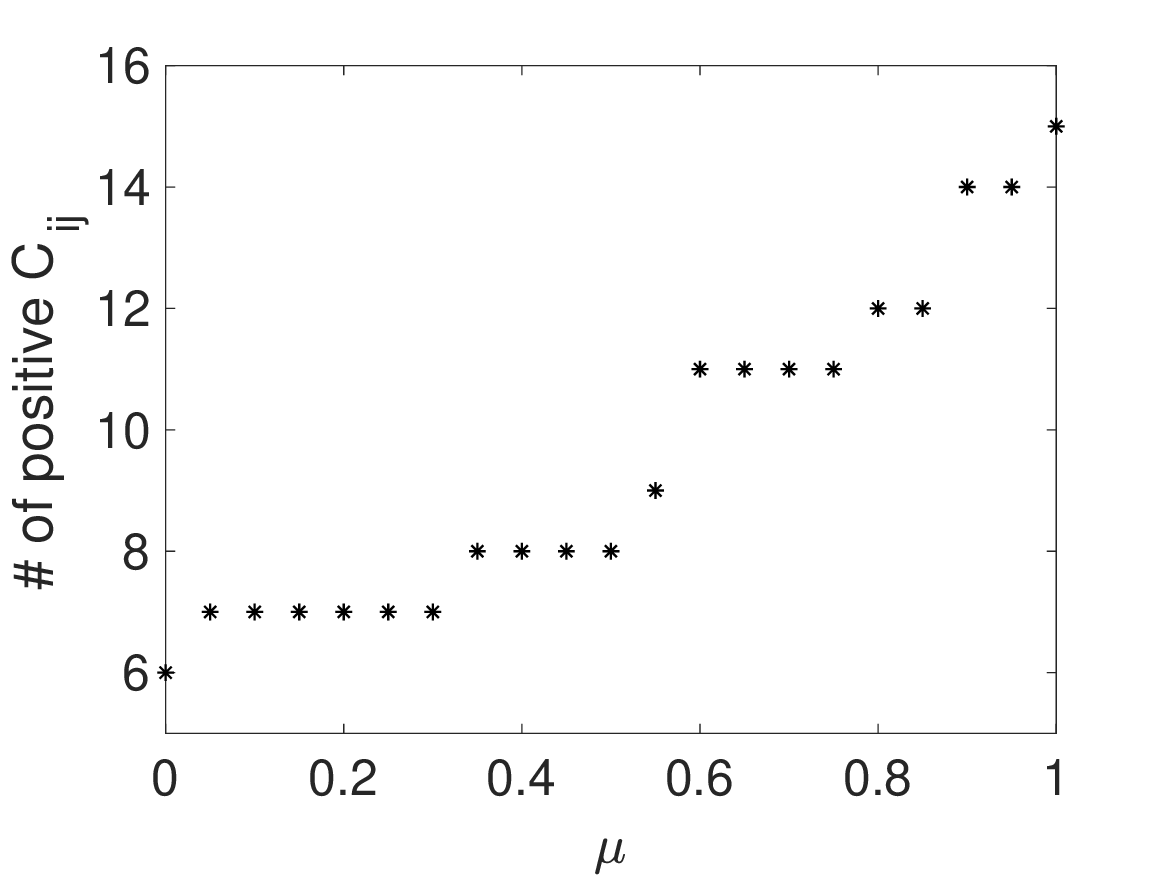} \\
 \caption{
 Results of minimization of the functional $\F=\F[C]$, given by \eqref{eq:F}, for the graph with $7$ vertices (Table \ref{Tab:1}), for $\mu\in[0,1]$. Top left panel: 
 Values of the energy functionals $\E=\E[C]$, given by \eqref{eq:energy}, and $\F=\F[C]$, for the minimizers $C\in\Cset$. 
Top right panel: Values of the kinetic (star) and metabolic (circle) energy of the minimizers, as defined in \eqref{eq:Ekinmet}.
Bottom left panel: Smallest three nonzero eigenvalues of the matrix Laplacian $\Lap[C]$. The second eigenvalue is the Fiedler number of the minimizer $C\in\Cset$.
Bottom right panel:
Number of nonzero elements of $C\in\Cset$, i.e., number of active edges of the minimizing graph.
\label{fig:N7En}}
\end{figure}

Statistical properties of the optimal graphs in dependence on the value of $\mu\in[0,1]$
are plotted in Fig. \ref{fig:N7En}.
In the top left panel, we plot the values of the energy $\E=\E[C]$ given by \eqref{eq:energy} and the modified energy $\F=\F[C]$ given by \eqref{eq:F}.
We observe that the value of $\E[C]$ is increasing with increasing $\mu\geq 0$. This can be seen as the extra energy expenditure for securing the robustness of the network.
In the top right panel we plot the kinetic (star-shaped markers) and metabolic (circular markers)
energies of the optimal networks, defined in \eqref{eq:Ekinmet}. The kinetic energy is decreasing with increasing $\mu$, which is due to the fact that less pumping power is necessary if the network consists of more edges, or edges with higher conductivities. However, this means that the ``material expense" is higher, i.e., the metabolic energy is increasing.
In the bottom left panel we plot the smallest three nonzero eigenvalues (i.e., the second, third and fourth smallest) of the matrix Laplacian $\Lap[C]$. This is to understand the multiplicity of the Fiedler number $\f=\f[C]$. We observe that the Fiedler number is simple for $\mu\lesssim 0.6$, turning to double for $\mu\gtrsim 0.6$.
Finally, in the bottom left panel of Fig. \ref{fig:N7En} the number of positive $C_{ij}$ (i.e., the number of edges present in the graph) is plotted. We again notice the monotonicity, i.e., increasing the value of $\mu$ indeed leads to the graph becoming denser.

\subsection{Leaf example}\label{subsec:leaf}
\begin{figure}[h] \centering
 \includegraphics[width=.5\textwidth]{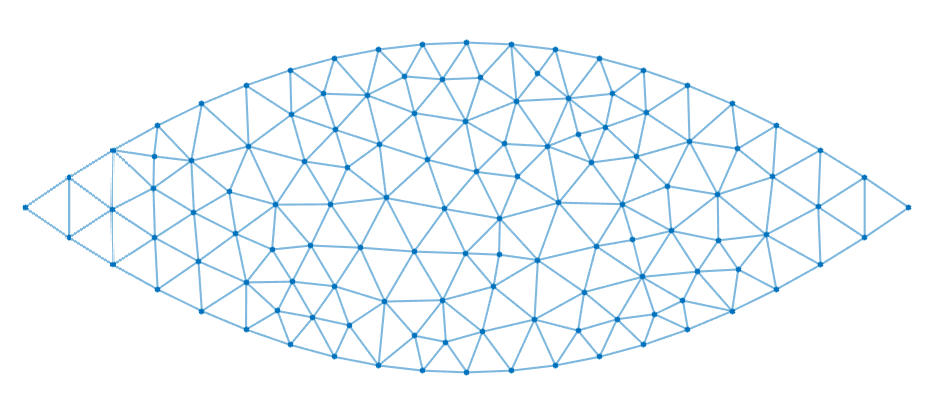}
 \caption{The planar graph modeling a leaf, with 122 nodes and 323 edges. \label{fig:leaf0}}
\end{figure}

Inspired by the possible application of the model to simulate leaf venation patterns,
we generated a planar graph in the form of a triangulation of leaf-shaped domain
with $|V|=122$ nodes and $323$ edges, see Fig. \ref{fig:leaf0}.
We prescribed a single source, $S_i=1$, for the left-most node
(``stem" of the leaf), while $S_j = -(|V|-1)^{-1}$ for all other nodes.

We again applied the projected subgradient method \eqref{eq:subgrad_step}--\eqref{eq:proj_step}
for minimization of the functional \eqref{eq:F} with $\gamma:=1$ and $\nu:=1$. We chose $\tau_0:=10^{-1}$ and $K:=10^7$, and $\mu\in [0,5]$. Although Lemma \ref{lem:boundbelow} guarantees boundedness from below of $\F=\F[C]$ only for $\mu\leq 1$, this condition is sufficient, but not necessary. As $\F=\F[C]$ is obviously bounded from below on every compact subset of $\Cset$, the projected subgradient method would indicate a possible unboundedness from below by divergence of the sequence of iterates. This, however, did not happen for $\mu\in [0,5]$.

\begin{figure}[h]\centering
 \includegraphics[width=.32\textwidth]{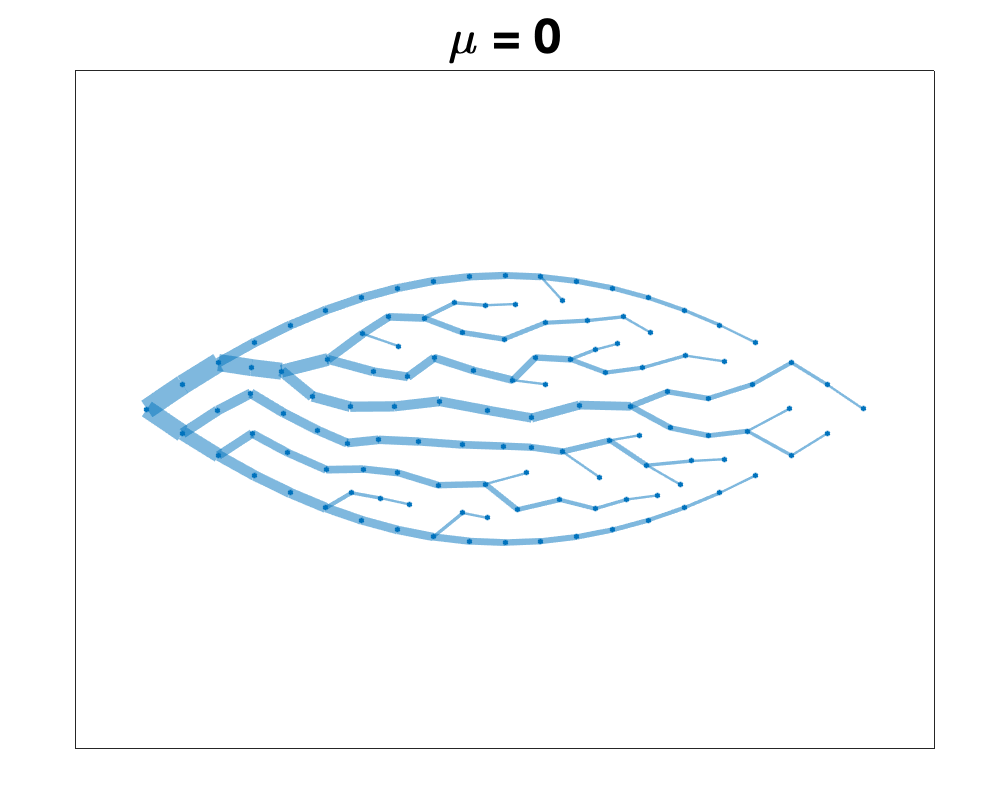}
 \includegraphics[width=.32\textwidth]{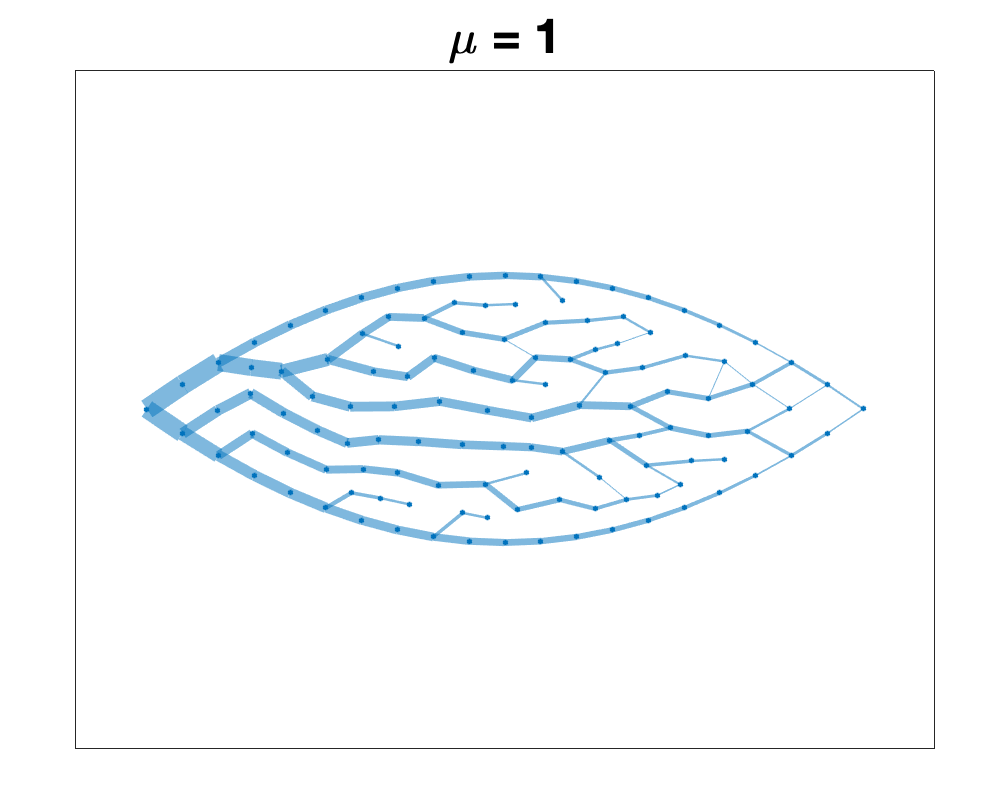}
 \includegraphics[width=.32\textwidth]{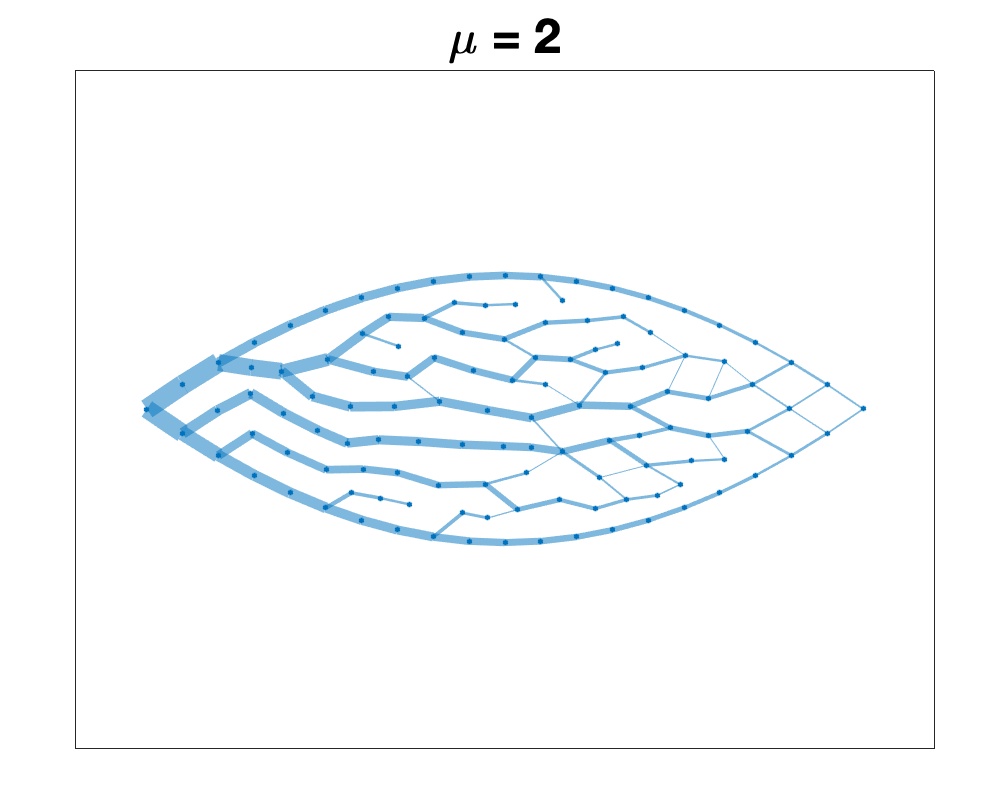} \\
 \includegraphics[width=.32\textwidth]{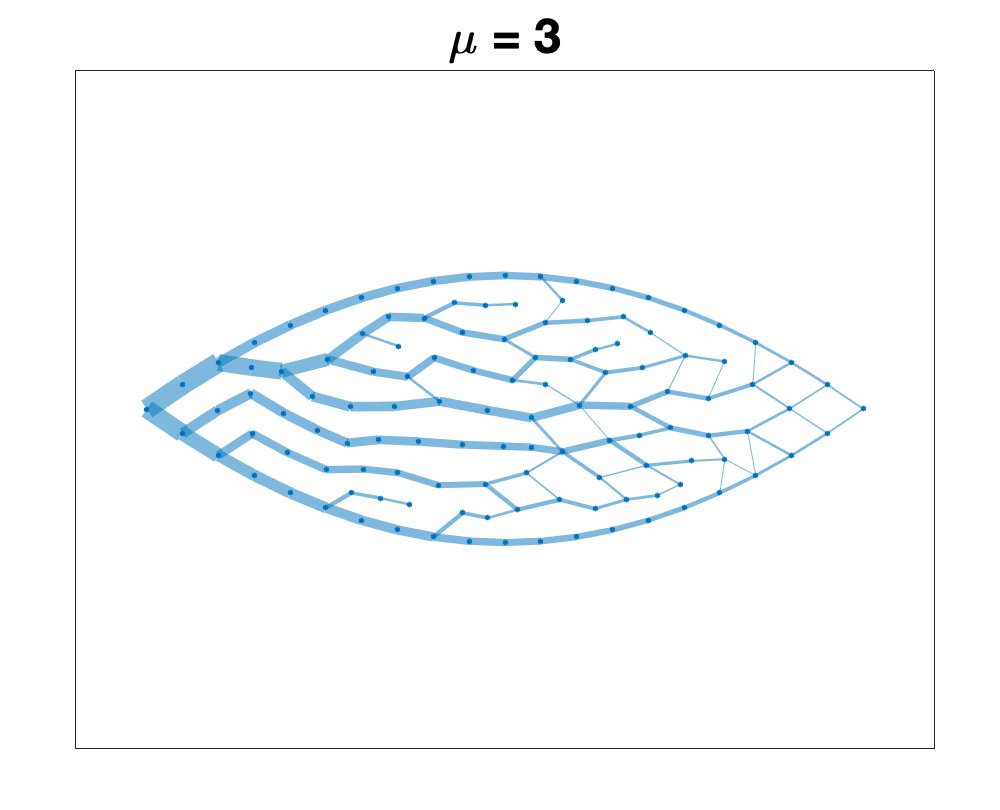}
 \includegraphics[width=.32\textwidth]{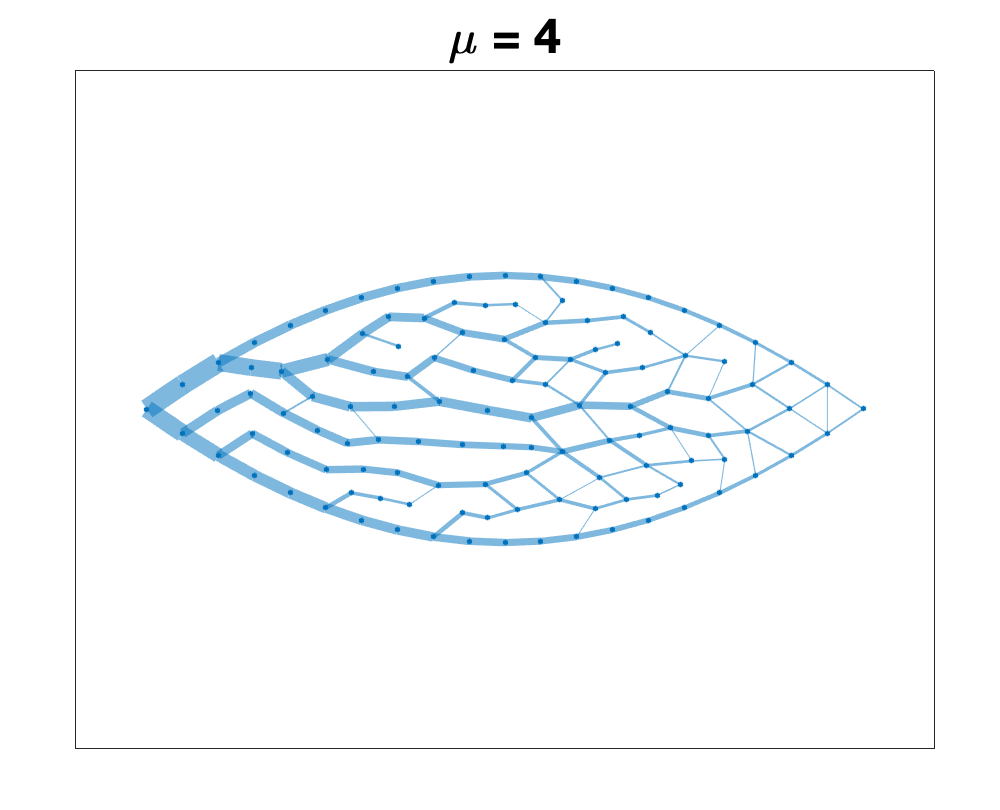}
 \includegraphics[width=.32\textwidth]{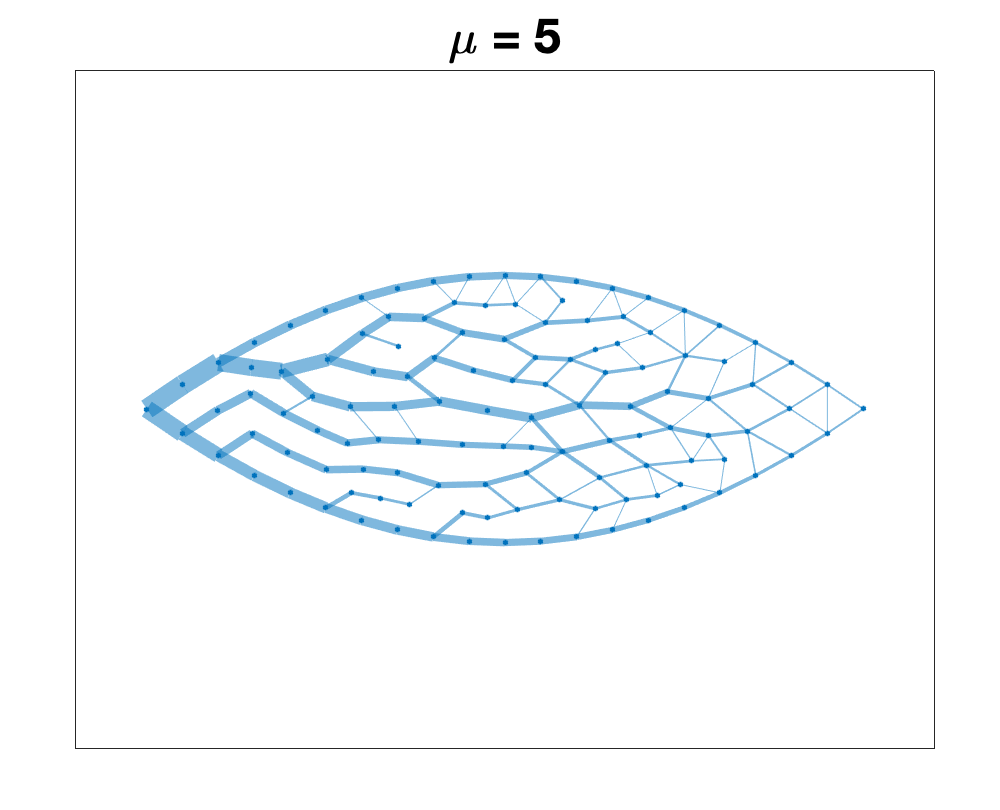} 
 \caption{Results of minimization of the functional $\F=\F[C]$, given by \eqref{eq:F}, for the leaf-shaped graph (Fig. \ref{fig:leaf0}) and $\mu\in\{0, 1, 2, 3, 4, 5\}$.
 The thickness of the line segments is proportional to the square root of the conductivity $C_{ij}\geq 0$ of the corresponding edge. Edges with $C_{ij}=0$ are excluded from the plot. \label{fig:leaf}}
\end{figure}

The optimal graphs found for $\mu\in\{0, 1, 2, 3, 4, 5\}$ are plotted in the correspondingly labeled panels of Fig. \ref{fig:leaf}.
Again, the thickness of the line segments is proportional to the square root of the conductivity of the corresponding edge. Edges with $C_{ij}=0$
are excluded from the plot.

\begin{figure}[h] \centering
 \includegraphics[width=.35\textwidth]{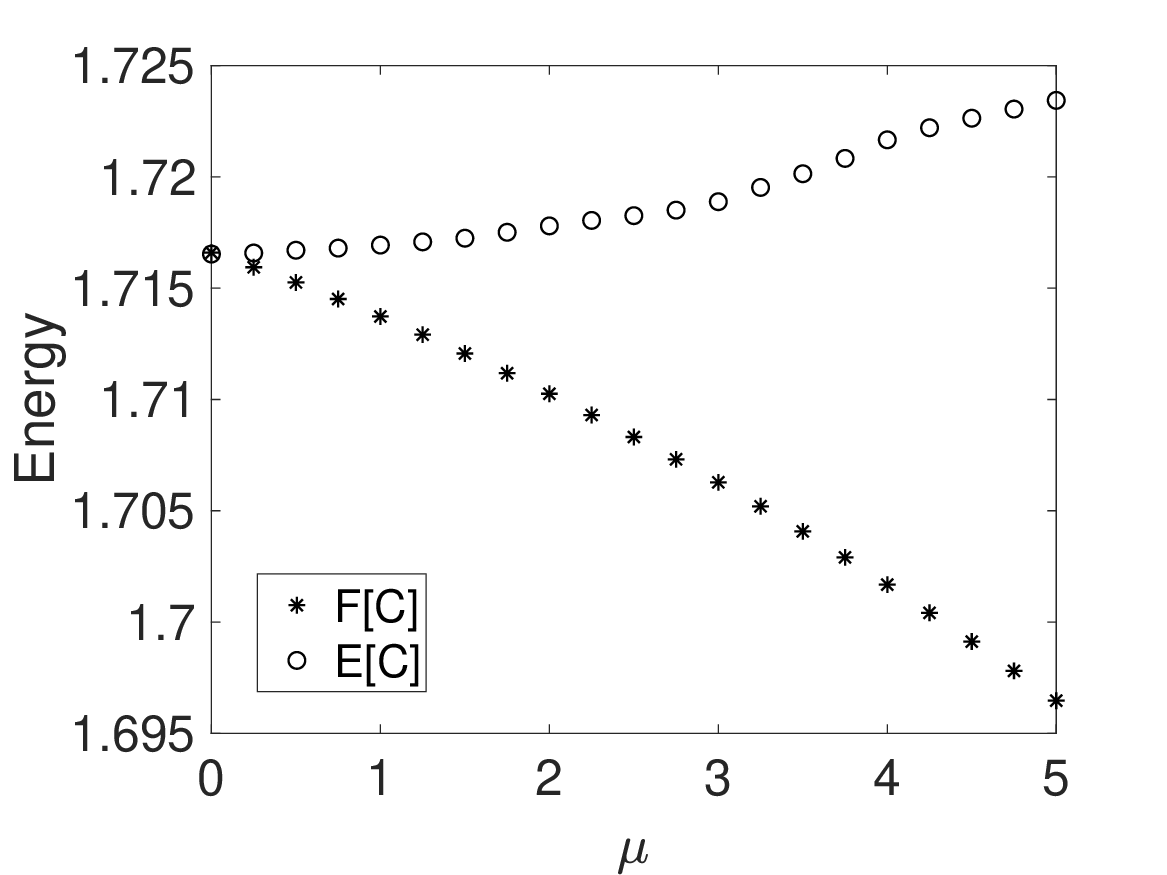}  $\qquad$
 \includegraphics[width=.35\textwidth]{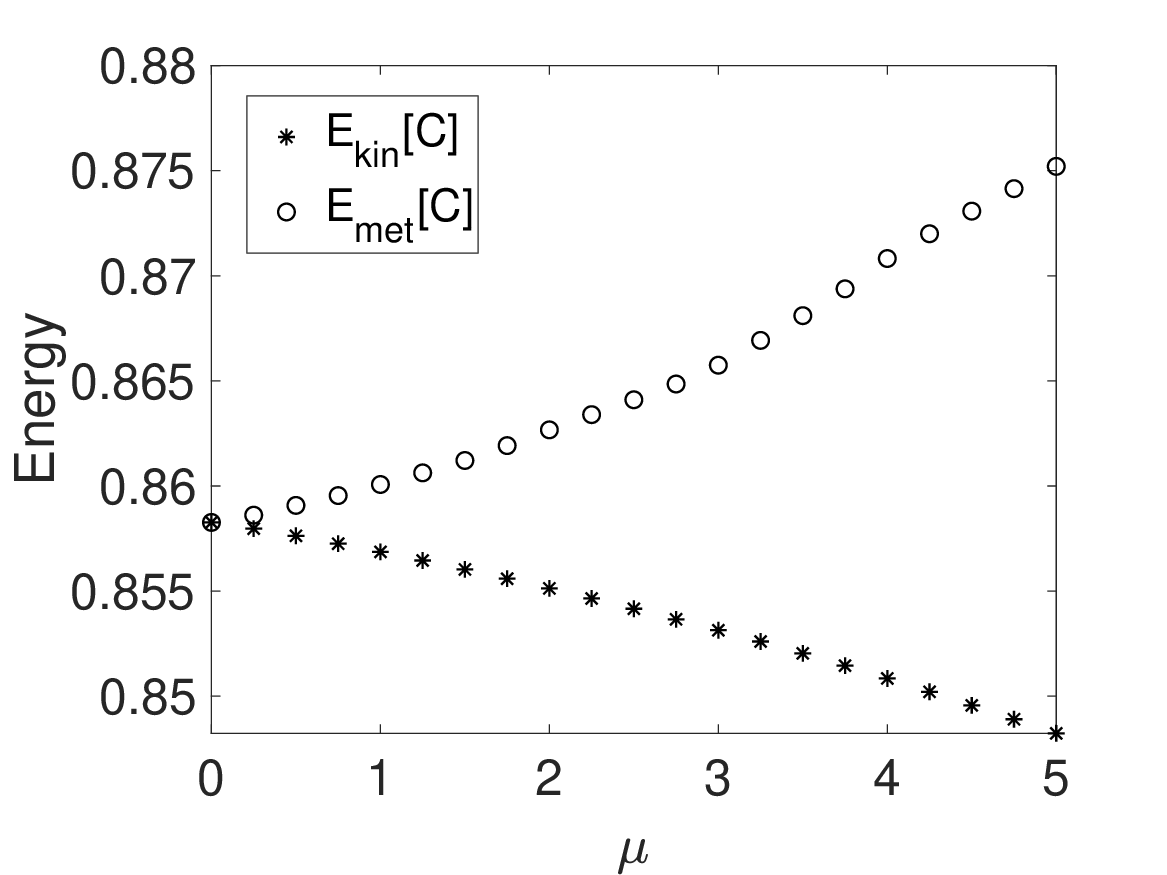} \\
 \includegraphics[width=.35\textwidth]{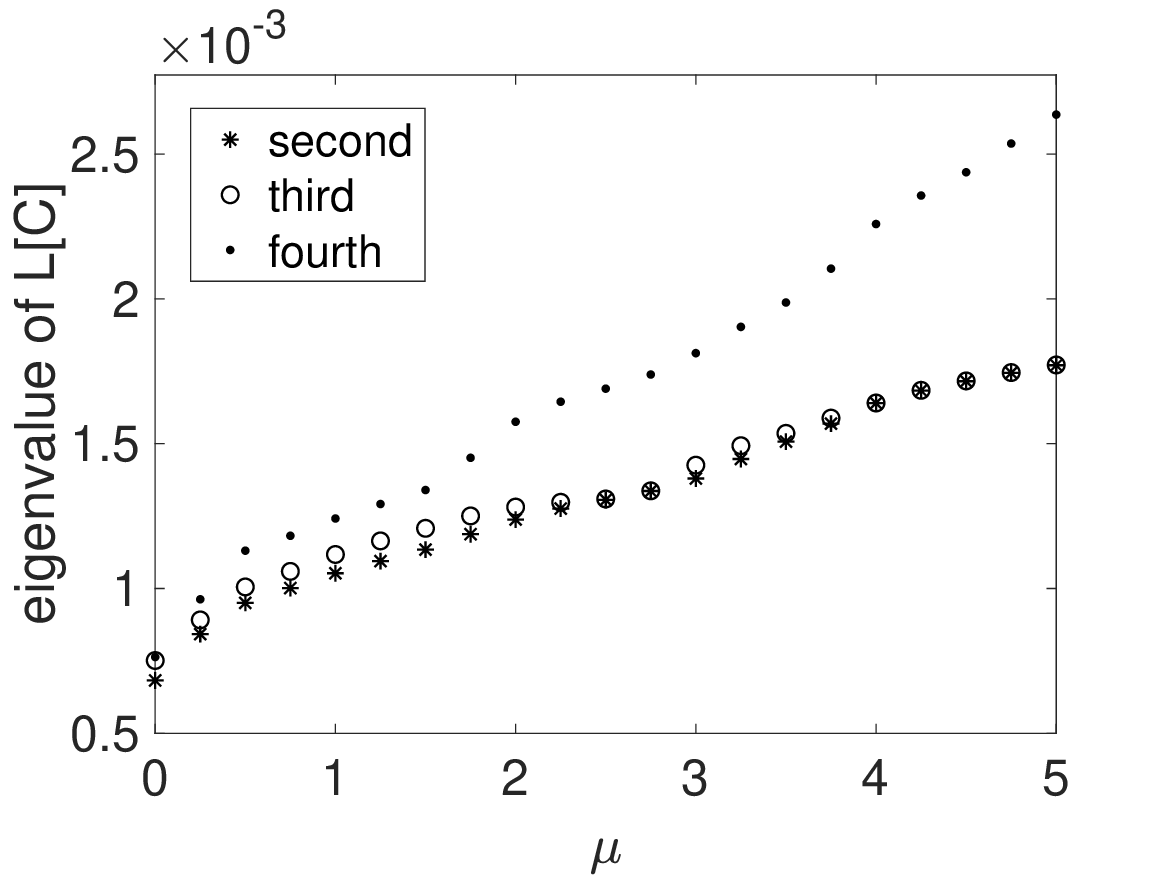}  $\qquad$
 \includegraphics[width=.35\textwidth]{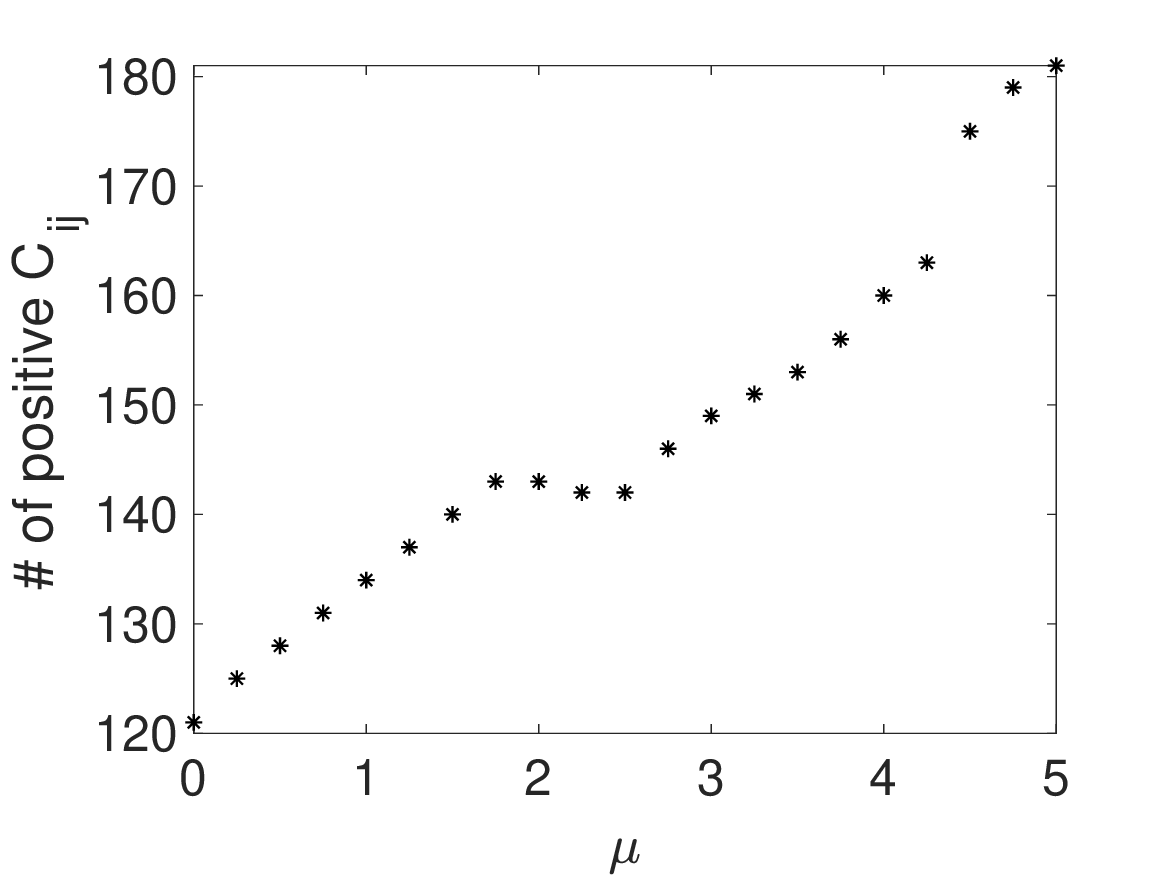} \\
 \caption{Results of minimization of the functional $\F=\F[C]$, given by \eqref{eq:F}, for the leaf-shaped graph (Fig. \ref{fig:leaf0}), for $\mu\in[0,5]$. Top left panel: 
 Values of the energy functionals $\E=\E[C]$, given by \eqref{eq:energy}, and $\F=\F[C]$, for the minimizers $C\in\Cset$. 
Top right panel: Values of the kinetic (star) and metabolic (circle) energy of the minimizers, as defined in \eqref{eq:Ekinmet}.
Bottom left panel: Smallest three nonzero eigenvalues of the matrix Laplacian $\Lap[C]$. The second eigenvalue is the Fiedler number of the minimizer $C\in\Cset$.
Bottom right panel:
Number of nonzero elements of $C\in\Cset$, i.e., number of active edges of the minimizing graph.
\label{fig:leafEn}}
\end{figure}

Statistical properties of the graphs in dependence on the value of $\mu\in [0,5]$ are plotted in Fig. \ref{fig:leafEn}.
Here, in the bottom left panel it is interesting to observe that the Fiedler number seems to be a double eigenvalue for $\mu\in\{2.5, 2.75\}$, and the again for $\mu\gtrsim 0.4$. In the bottom right panel we see that the number of active edges is no more monotonically increasing with $\mu$. Indeed, around the value of $\mu\simeq 2.5$ the number of nonzero elements of $\widetilde C$ drops slightly. We hypothesize that this ``anomaly" may be related to the multiplicity of the Fiedler number observed in the left bottom panel around the same value of $\mu$. However, we are not able to provide a rigorous explanation.

\section{Appendix}\label{sec:app}

\subsection{Some properties of the energy functional \eqref{eq:energy}}

Here we prove some mathematical properties of the functional $\E=\E[C]$ and the set ${\mathcal C}$, defined by \eqref{eq:energy} and \eqref{eq:setC}, respectively.

\begin{lemma}\label{lem:C1}
For any fixed $S\in\R^{|V|}_0$,
the set
\begin{equation}\label{eq:Efinite}
\{C\in \Cset:\E[C]<+\infty\}
\end{equation}
is an open convex subset of $\Cset$.
\end{lemma}

\begin{proof}
We start by showing that \eqref{eq:Efinite} is convex. For that sake, we take $C^1,C^2\in \Cset$ with $\E[C^1]<+\infty$ and $\E[C^2]<+\infty$.
Let $0<\alpha<1$ and put $C:=\alpha C^1+(1-\alpha)C^2.$ We want to show that also $\E[C]<+\infty$, i.e., that \eqref{eq:kirchhoff}
is solvable with $C$.

Let $x\in\R^{|V|}$ be fixed. Then
\begin{equation}\label{eq:Lapl1}
x^T\Lap[C]x=\alpha x^T\Lap[C^1]x+(1-\alpha)x^T\Lap[C^2]x.
\end{equation}
Since all the three Laplacian matrices appearing in the above identity are positive semi-definite, we observe that the left-hand side of \eqref{eq:Lapl1}
is zero if, and only if, both terms on the right-hand side of \eqref{eq:Lapl1} vanish. Therefore, the null space of $\Lap[C]$ is the intersection
of the null spaces of $\Lap[C^1]$ and $\Lap[C^2]$ and the range space of $\Lap[C]$ is the sum of range spaces of $\Lap[C^1]$ and $\Lap[C^2].$ We conclude that $S$ is
in the range space of $\Lap[C]$ and \eqref{eq:kirchhoff} is solvable for $C$.

Next, we show that \eqref{eq:Efinite} is open in $\Cset$. Let $C\in\Cset$ with $\E[C]<+\infty$ and let $\Gamma_1,\dots,\Gamma_n\subset V$ be the connected components of $C$.
Then \eqref{eq:kirchhoff} is solvable if, and only if, $\sum_{j\in\Gamma_k}S_j=0$ for every $1\le k\le n.$
Let us define
\[
   \mathcal{U}[C] :=
   \left\{ C'\in\Cset; \, C_{ij} > 0 \implies C'_{ij} > 0 \mbox{ for all } (i,j)\in\Eset \right\}.
\]
Obviously, any $C'\in \mathcal{U}[C]$ has either the same
connected components as $C$, or fewer but larger components, stemming from establishing new connections among $\Gamma_1,\dots,\Gamma_n\subset V$.
In both cases \eqref{eq:kirchhoff} is solvable with $C'$ and, therefore, $\E[C']$ is also finite. Consequently, $\mathcal{U}[C]$ is an open neighborhood of $C$ in the set $\{C\in \Cset:\E[C]<+\infty\}$.
\end{proof}

In the sequel we shall denote, for any matrix $A\in\R^{|V|\times |V|}$ and any vector $x\in\R^{|V|}$,
\[
   \Norm{A}_\infty := \max_{i, j \in \Vset} |A_{ij}|, \qquad
   |x|_\infty := \max_{i \in \Vset} |x_{i}|.
\]
Moreover, we recall that $\R^{|V|}_0$ was defined in \eqref{eq:RV0} as the set of vectors $x\in\R^{|V|}$
with vanishing sum.

Next we study the continuity of the pressures $P=P[C]$, the fluxes $Q=Q[C]$ and the energy $\E=\E[C]$ as functions of $C$. 
We restrict ourselves to the case when the Kirchhoff law \eqref{eq:kirchhoff} has a unique solution $P=P[C]\in\R^{|V|}_0$.
This happens if, and only if, the matrix Laplacian $\Lap[C]$ of $C$ is nonsingular on $\R_0^{|V|}$,
which in turn holds if, and only if, $C$ represents a weighted connected graph.


\begin{example}\label{ex:P}
If $C$ represents an unconnected graph, then \eqref{eq:kirchhoff} might admit no solution 
at all or the solution might not be unique in $\R^{|V|}_0$.
In both cases we might interpret $P=P[C]$ as a set-valued function of $C$, cf. \cite{AF_set}.
Unfortunately, the following example shows that $P$ does not need to be continuous.

Let $G=(V,E)$ be a complete graph on four vertices $V=\{1,2,3,4\}$.
Let $L_{ij}=1$ for $i\not=j$, and let $S=(1,-1,1,-1)$. Furthermore, set
\[C^0:=\left(\begin{matrix}0&1&0&0\\1&0&0&0\\0&0&0&1\\0&0&1&0\end{matrix}\right).\]
Then the graph represented by $C^0$ is disconnected with two components $U_1=\{1,2\}$ and $U_2=\{3,4\}$. Since  
$S_1+S_2=S_3+S_4=0$, the equation
$\Lap[C^0]P=S$ has non-unique solutions, namely every $P\in\{(a,a-1,b,b-1), a,b\in\R\}$.
Even when we restrict $P$ to be an element of $\R_0^4$, we still have the non-unique solutions $P=(a,a-1,1-a,-a)$
with $a\in\R$. 

For $t>0$, we define $C^t$ by adding to $C^0$ the edge $(2,3)$ with weight $t$. Namely,
we define $C^t_{2,3}=C^t_{3,2}=t$ and $C^t_{ij}=C^0_{ij}$ for $\{i,j\}\not=\{2,3\}.$
Then $\Lap[C^t]P=S$ for every $P=(c+1,c,c,c-1)$, $c\in\R$ and in $\R_0^4$, the solution is unique, namely 
$P=(1,0,0,-1)$. 

On the other hand, if $t<0$ we define $C^t$ by adding an edge $(1,4)$ with a weight $-t>0$ to $C^0$.
Then $\Lap[C^t]P=S$ for every
$P=(c,c-1,c+1,c)$, $c\in\R$. Again, in $\R_0^4$, the solution is unique, namely $P=(0,-1,1,0)$.
Using the notions of set-valued analysis \cite{AF_set}, we observe that
the mapping $t\to \{P\in\R_0^4:\Lap[C^t]P=0\}$ is upper semi-continuous
but not lower semi-continuous in $t=0$.
\end{example}

We formulate the next result as a general statement about matrix Laplacian $\Lap[K]$ of a symmetric matrix $K$ with non-negative
entries. Later on, we apply this result to the matrix $K_{ij}=C_{ij}/L_{ij}$, $(i,j)\in E$.

\begin{lemma} \label{lemma:pert}
Let $S\in \R^{|V|}_0$, let $K^\ast\in \R_+^{|V|\times|V|}$ be a symmetric matrix with non-negative entries representing a connected weighted undirected graph
and let $P^\ast\in \R^{|V|}_0$ be the unique solution of the linear system
\( \label{LastPastSast}
   \Lap(K^\ast) P^\ast = S.
\)
Then there exists a small neighbourhood of $K^\ast$ (relative to $\R_+^{|V|\times|V|}$) on which \eqref{LastPastSast} is uniquely solvable
and defines $P=P(K)$ as a differentiable, Lipschitz-continuous function of $K$. Moreover, 
the fluxes $Q_{ij}=K_{ij} (P_j - P_i)$ are Lipschitz-continuous functions of $K$ as well.
\end{lemma}

\begin{remark}
To be more specific, Lemma \ref{lemma:pert} ensures that to a given $K^\ast\in \R_+^{|V|\times|V|}$
there exists $\rho>0$, which in general depends on $K^\ast$, such that if
$K\in \R_+^{|V|\times|V|}$ verifies
$\Norm{K - K^\ast}_\infty \leq \rho,$
then there exists a unique solution $P\in \R^{|V|}_0$ of the linear system
$\Lap(K) P = S.$ 
Moreover, there exists a constant $c$, independent of $\rho$, such that the fluxes
\[
   Q_{ij} := K_{ij} (P_j - P_i), \qquad Q_{ij}^\ast := K_{ij}^\ast (P_j^\ast - P_i^\ast)
\]
satisfy
\(   \label{QQast}
   \Norm{Q-Q^\ast}_\infty \leq c \Norm{K - K^\ast}_\infty.
\)
\end{remark}

\begin{proof}
Since $\Lap(K^\ast)$ is a nonsingular operator on the space $\R^{|V|}_0$, a classical perturbation
result implies unique solvability of
\(  \label{LPS}
   \Lap(K) P = S
\)
on $\R^{|V|}_0$ for $K\in\R_+^{|V|\times|V|}$ with
\(   \label{CCast}
    \Norm{K - K^\ast}_\infty \leq \rho
\)
and $\rho$ small enough. Let us denote $\delta K := K-K^\ast$ and $\delta \Lap := \Lap(K)-\Lap(K^\ast)$.
Note that $\delta\Lap$ depends linearly on $\delta K$, $\Norm{\delta \Lap}_\infty \leq (|V|-1)\Norm{\delta K}_\infty$ and that $\delta\Lap$ maps $\R_0^{|V|}$ into itself.

In particular, denoting $\delta P := P-P^\ast$, we have
\[
   (\Lap(K^\ast) + \delta \Lap)(P^\ast + \delta P) = S.
\]
Using \eqref{LastPastSast}, applying $\Lap(K^\ast)^{-1}$ to both sides of this equation and adding $P^\ast$ gives
\[
P^\ast+(\Lap(K^\ast))^{-1} \delta \Lap P^\ast + \delta P+(\Lap(K^\ast))^{-1} \delta \Lap\delta P=P^\ast,
\]
which can be reformulated as
\(  \label{dP}
   \delta P = \left[ \left( I + (\Lap(K^\ast))^{-1} \delta \Lap \right)^{-1} - I \right] P^\ast.
\)
Expanding the expression $\left( I + (\Lap(K^\ast))^{-1} \delta \Lap \right)^{-1}$ into the von Neumann series, we have
\begin{equation}\label{eq:Neumann}
\delta P=\sum_{k=1}^\infty [-(\Lap(K^\ast))^{-1} \delta \Lap]^k P^\ast.
\end{equation}
As $\delta\Lap$ depends linearly on $\delta K$, it follows that $P=P(K)$ is differentiable at $K^\ast.$

In the next step, we denote by $\Norm{\cdot}$ the operator norm induced by the vector norm $|\cdot|_\infty$ on the space $\R^{|V|}_0$.
By norm equivalence, there exists a constant $c>0$ such that
$\Norm{\cdot} \leq c \Norm{\cdot}_\infty$ and we recall that $\Norm{\delta \Lap}_\infty \leq (|V|-1)\Norm{\delta K}_\infty$.
Combining this with \eqref{eq:Neumann} we obtain
\begin{align}  \label{deltaP}
   |\delta P|_\infty &\leq \sum_{k=1}^\infty \|(\Lap(K^\ast))^{-1} \delta \Lap\|^k |P^\ast|_\infty
\leq |P^\ast|_\infty\cdot \sum_{k=1}^\infty \|(\Lap(K^\ast))^{-1}\|^k\cdot \|\delta \Lap\|^k \\
&\leq  |P^\ast|_\infty\cdot \sum_{k=1}^\infty [c(|V|-1)\|(\Lap(K^\ast))^{-1}\|]^k\|\delta K\|_\infty^k
\notag= |P^\ast|_\infty\cdot\frac{\Lambda \Norm{\delta K}_\infty}{1 - \Lambda \Norm{\delta K}_\infty},
\end{align}
where we denoted $\Lambda := c \Norm{(\Lap(K^\ast))^{-1}} (N-1)$ and assumed that
$\Norm{\delta K}_\infty < \Lambda^{-1}$.

Now, denoting $(\Delta P)_{ij} := P_j - P_i$ and analogously for $(\Delta P^\ast)_{ij}$, we have for every $i, j \in \Vset$,
\[
   \bigl| Q_{ij} - Q_{ij}^\ast \bigr| &=& \bigl| K_{ij}(\Delta P)_{ij} - K_{ij}^\ast (\Delta P^\ast)_{ij} \bigr| \\
     &\leq& \bigl| (K_{ij} - K_{ij}^\ast)((\Delta P)_{ij} - (\Delta P^\ast)_{ij}) \bigr|
         +  \bigl| (K_{ij} - K_{ij}^\ast) (\Delta P^\ast)_{ij}) \bigr| + \bigl| K_{ij}^\ast ((\Delta P)_{ij} - (\Delta P^\ast)_{ij}) \bigr| \\
   &\leq& 2 \Norm{\delta K}_\infty |\delta P|_\infty + \Norm{\delta K}_\infty \Norm{\Delta P^\ast}_\infty + 2 \Norm{K^\ast}_\infty |\delta P|_\infty,
\]
where we used the estimate $| (\Delta P)_{ij} - (\Delta P^\ast)_{ij} | \leq 2 |\delta P|_\infty$.
Therefore, using \eqref{CCast} and \eqref{deltaP}, we estimate
\[
   \bigl| Q_{ij} - Q_{ij}^\ast \bigr| \leq
      2 \left(\Norm{\delta K}_\infty + \Norm{K^\ast}_\infty \right) \frac{\Lambda \Norm{\delta K}_\infty}{1 - \Lambda \Norm{\delta K}_\infty} + \Norm{\delta K}_\infty \Norm{\Delta P^\ast}_\infty.
\]
Finally, we choose $\rho := 1/(2\Lambda)$, so that for $\Norm{\delta K}_\infty \leq \rho$ we have
\[
    \frac{1}{1 - \Lambda \Norm{\delta K}_\infty} \leq 2, \qquad \Lambda \Norm{\delta K}_\infty^2 \leq \frac12 \Norm{\delta K}_\infty,
\]
which gives
\[
   \bigl| Q_{ij} - Q_{ij}^\ast \bigr| \leq \left( 2 + 4 \Lambda \Norm{K^\ast}_\infty + \Norm{\Delta P^\ast}_\infty \right) \Norm{\delta K}_\infty
\]
and an obvious choice of the constant $c$ concludes \eqref{QQast}.
\end{proof}

\begin{corollary} 
The functional $\E=\E[C]$ is continuous on $\Cset$ and
totally differentiable on the set $\{C\in\Cset: C\ \text{corresponds to a connected graph}\}\subset\Cset$.
\end{corollary}

\begin{proof}
By Lemma \ref{lem:kin_convex}, $\E$ is continuous in $C$ whenever $\E[C]<+\infty$.
It is therefore enough to show that $\E(C)$ is large on a small neighbourhood of a given $C^\ast$ with $\E(C^\ast)=+\infty.$
In that case, \eqref{eq:kirchhoff}
is not solvable and, therefore, $C^\ast$ represents a disconnected graph. Let us assume for simplicity, that it has
only two connected components $\Gamma_1, \Gamma_2 \subset \Vset$ and define $A:=\sum_{j\in \Gamma_1}S_j=-\sum_{j\in \Gamma_2}S_j >0$.
If \eqref{eq:kirchhoff} is solvable for $C\in \Cset$ with $\|C-C^\ast\|_\infty<\eps$, then there must be
some new edges between $\Gamma_1$ and $\Gamma_2$ with weights bounded by $\eps$,
that transfer the mass $A$ from $\Gamma_1$ to $\Gamma_2$.
Therefore, $\Ekin(C)\ge \min_{(i,j)\in E} A^2 L_{i,j}/(\eps |V|^2)$, which tends to infinity as $\eps\to 0$.
We thus conclude that $\E=\E[C]$ is continuous on $\Cset$.

The statement about total differentiability of $\E=\E[C]$ follows directly from Lemma \ref{lemma:pert}.
\end{proof}

\subsection{Derivative of the Fiedler number}\label{subsec:derivative}

Here we prove the following result.

\begin{lemma}\label{lem:DFiedler}
Assume that the Fiedler number $\f[C]$ of the conductivity matrix $C\in\Cset$ is a simple eigenvalue of the matrix Laplacian $\Lap[C]$, given in Definition \ref{def:Lapl}.
Then
\(   \label{eq:DFiedler}
    \part{\f[C]}{C_{ij}} = (v_i - v_j)^2,
\)
where $v$ is the normalized eigenvector of $\Lap[C]$ corresponding to the eigenvalue $\f[C]$.
\end{lemma}

We start with an auxiliary result about the derivative of a simple eigenvalue of a symmetric matrix with respect to the matrix elements. This is a well-known result going back at least to Rellich \cite{Rellich}, see also \cite{Kato, Stewart-Sun}. However, we offer its simple proof here for the sake of the reader.

\begin{lemma} \label{lem:Dlambda}
Let $\lambda\in\R$ be a simple eigenvalue of a symmetric matrix $\Lap\in\R^{d\times d}$, $d>0$, with normalized eigenvector $v\in\R^d$.
Then
\( \label{eq:Dlambda}
   \part{\lambda}{\Lap_{ij}}  = \begin{cases}
       v_i^2 \qquad\mbox{for } i = j, \\
        2 v_i v_j \qquad\mbox{for } i \neq j,
        \end{cases}
\)
where we take the variation in the set of symmetric matrices, i.e.,
in both the elements $\Lap_{ij}$ and $\Lap_{ji}$.
\end{lemma}

\begin{proof}
We differentiate the identity $\Lap v = \lambda v$ with respect to $\Lap_{ij}$,
\[
   \sum_{m=1}^d \part{\Lap_{km}}{\Lap_{ij}} v_m + \sum_{m=1}^d \Lap_{km} \part{v_m}{\Lap_{ij}} = \part{\lambda}{\Lap_{ij}} v_k + \lambda \part{v_k}{\Lap_{ij}}
\]
for all $k=1,\dots,d$.
Multiplication by $v_k$ and summation in $k$ gives
\[
   \sum_{k=1}^d \sum_{m=1}^d \part{\Lap_{km}}{\Lap_{ij}} v_m v_k + \sum_{k=1}^d \sum_{m=1}^d \Lap_{km} \part{v_m}{\Lap_{ij}} v_k = \part{\lambda}{\Lap_{ij}} +  \lambda  \sum_{k=1}^d \part{v_k}{\Lap_{ij}} v_k,
\]
where we used the normalization $\Norm{v}^2=1$ for the first term of the right-hand side.
Moreover, differentiation of this identity with respect to $\Lap_{ij}$ gives $\sum_{k=1}^d \part{v_k}{\Lap_{ij}} v_k = 0$, so that
\[
   \sum_{k=1}^d \sum_{m=1}^d \part{\Lap_{km}}{\Lap_{ij}} v_m v_k + \sum_{k=1}^d \sum_{m=1}^d \Lap_{km} \part{v_m}{\Lap_{ij}} v_k = \part{\lambda}{\Lap_{ij}}
\]
Due to the symmetry of $L$, we have
\[
   \sum_{k=1}^d \sum_{m=1}^d \Lap_{km} \part{v_m}{\Lap_{ij}} v_k = \lambda \sum_{m=1}^d \part{v_m}{\Lap_{ij}} v_m = 0
\]
and, finally,
\[
   \sum_{k=1}^d \sum_{m=1}^d \part{\Lap_{km}}{\Lap_{ij}} v_m v_k = \begin{cases}
        v_i^2 \qquad\mbox{if } i=j, \\
      2 v_i v_j \qquad\mbox{if } i\neq j.
      \end{cases}
\]
\end{proof}

Now we easily evaluate the derivative of the Fiedler number $\f=\f[C]$, given it is a simple eigenvalue, with respect to the elements of the conductivity matrix $C_{ij}$, which amounts to the proof of Lemma \ref{lem:DFiedler}.

\begin{proof}
Let us fix any $i\neq j$. We have
\[
    \part{\f[C]}{C_{ij}} = \sum_{k=1}^d\sum_{m\neq k} \part{\f[C]}{\Lap_{km}} \part{\Lap_{km}}{C_{ij}} + \sum_{k=1}^d \part{\f[C]}{\Lap_{kk}} \part{\Lap_{kk}}{C_{ij}},
\]
where here and in the sequel we write $\Lap$ for $\Lap[C]$.
By Definition \ref{def:Lapl},
\[
   \part{\Lap_{kk}}{C_{ij}} = \sum_{m=1}^d \part{C_{km}}{C_{ij}} = \delta_{k,i} + \delta_{k,j},
\]
so that \eqref{eq:Dlambda} gives
\[
   \sum_{k=1}^d \part{\f[C]}{\Lap_{kk}} \part{\Lap_{kk}}{C_{ij}} = v_i^2 + v_j^2.
\]
Moreover, for $k\neq m$ we have $\part{\Lap_{km}}{C_{ij}} = -\delta_{k,i}\delta_{m,j}$, so that
\[
   \sum_{k=1}^d\sum_{m\neq k} \part{\f[C]}{\Lap_{km}} \part{\Lap_{km}}{C_{ij}} = - 2 v_i v_j,
\]
and \eqref{eq:DFiedler} follows.
\end{proof}


\end{document}